\documentclass[a4paper,11pt,onecolumn,twoside]{article}
\usepackage{mathrsfs}
\usepackage{mathrsfs}
\usepackage{amsfonts}
\usepackage{fancyhdr}
\usepackage{amsmath,amsfonts,amssymb,graphicx,amsthm}
\allowdisplaybreaks

\usepackage{amsmath,amsfonts,amssymb,amscd,amsthm,amsbsy,bm,latexsym}
\usepackage{gensymb}
\usepackage[pdftex,bookmarksnumbered=true,bookmarksopen=true,          
            colorlinks=true,pdfborder=001,citecolor=blue,
            linkcolor=red,anchorcolor=green,urlcolor=blue]{hyperref}
\usepackage{mathrsfs}
\allowdisplaybreaks

\usepackage{color}

\begin{document}

\title{\bf Asymptotic stability of a composite wave  for the one-dimensional compressible micropolar fluid model without  viscosity}
\author{
{\bf   Liyun Zheng}\\
School of Mathematical Sciences,
Anhui University, Hefei 230601,  P. R. China\\[2mm]
{\bf Zhengzheng Chen}\thanks{Corresponding author. E-mail:
chenzzandu@163.com}\\
School of Mathematical Sciences,
Anhui University, Hefei 230601,  P. R. China\\[2mm]
{\bf   Sina Zhang}\\
School of Mathematical Sciences,
Anhui University, Hefei 230601,  P. R. China\\[2mm]
}

\date{}

\vskip 0.2cm

\maketitle

\vskip 0.2cm \arraycolsep1.5pt
\newtheorem{Lemma}{Lemma}[section]
\newtheorem{Theorem}{Theorem}[section]
\newtheorem{Definition}{Definition}[section]
\newtheorem{Proposition}{Proposition}[section]
\newtheorem{Remark}{Remark}[section]
\newtheorem{Corollary}{Corollary}[section]

\begin{abstract}
We are concerned with  the large time behavior of solutions to the Cauchy problem of the one-dimensional compressible  micropolar fluid model without  viscosity, where  the far-field states of the  initial data are prescribed  to be different. If the corresponding Riemann problem of the compressible Euler system admits a contact discontinuity and two rarefaction waves solutions, we show that for such a non-viscous model,  the combination of the viscous contact wave with two rarefaction waves is time-asymptotically stable provided that the strength of the composite wave and the initial perturbation are sufficiently small. The proof is given by an elementary $L^2$ energy method.

\bigbreak
\noindent

{\bf \normalsize Keywords}\,\, {Compressible micropolar fluid model;\,\,Viscous contact wave;\,\, Rarefaction waves;\\\,\,Without viscosity;\,\,Nonlinear Stability}
\bigbreak

\end{abstract}

\section{Introduction }
\setcounter{equation}{0}
The one-dimensional full compressible micropolar fluid model in the Lagrangian coordinates reads as:
\begin{eqnarray}\label{1.0}
\left\{\begin{array}{ll}
          v_t-u_x=0,\\[2mm]
          u_t+p_x=\displaystyle\left(\frac{\mu u_x}{v}\right)_x,\\[2mm]
          \displaystyle\left(e+\frac{u^{2}}{2}\right)_{t}+(pu)_{x}=\displaystyle\left(\frac{\kappa\theta_{x}}{v}\right)_{x}+\frac{\omega_{x}^{2}}{v}+v\omega^{2},\\[4mm]
          \omega_{t}=\displaystyle A\left[\left(\frac{w_{x}}{v}\right)_{x}-v\omega\right],
\end{array}\right. \,\,t>0,\, x\in \mathbb{R},
\end{eqnarray}
where the unknown functions are the specific volume  $v(t,x)>0$, the velocity $u(t, x)$, the microrotation velocity $w$, the absolute temperature $\theta(t, x)>0$, the internal energy $e(t, x)$ and the pressure $p(v,\theta)$ of the fluid, respectively, while  the constants $\mu, \kappa$ and $A$ denote the viscosity  coefficient, the heat conductivity coefficient and the microviscosity  coefficient, respectively.

The model of  micropolar fluid was first introduced  by Eringen \cite{A.C.Erigen-1966} in 1966. This model can be used  to describe the motions of a large variety of complex fluids consisting of dipole elements such as the suspensions, animal blood, liquid crystal, etc. For more physical background on this model, we  refer  to \cite{A.C.Erigen-1999, G.Lukaszewicz-1999}.  Note that, if the microstructure of the fluid is neglected, i.e.,$w=0$,  then system (\ref{1.0}) is reduced  to  the classical compressible Navier-Stokes system.

In this paper, we consider the the Cauchy problem of the system (\ref{1.0}) without viscosity, i.e.,
\begin{eqnarray}\label{1.1}
\left\{\begin{array}{ll}
          v_t-u_x=0,\\[2mm]
          u_t+p_x=0,\\[2mm]
          \displaystyle\left(e+\frac{u^{2}}{2}\right)_{t}+(pu)_{x}=\displaystyle\left(\frac{\kappa\theta_{x}}{v}\right)_{x}+\frac{\omega_{x}^{2}}{v}+v\omega^{2},\\[4mm]
          \omega_{t}=\displaystyle A\left[\left(\frac{w_{x}}{v}\right)_{x}-v\omega\right],
\end{array}\right. \,\,t>0,\, x\in \mathbb{R}
\end{eqnarray}
with the following initial and far field conditions:
\begin{eqnarray}\label{1.2}
\left\{\begin{array}{ll}
          (v,u,\theta,\omega)(x,0)=(v_{0},u_{0},\theta_{0},\omega_{0})(x),\quad x\in \mathbb{R},\\[2mm]
         (v,u,\theta,\omega)(\pm\infty,t)=(v_{\pm},u_{\pm},\theta_{\pm}, \omega_{\pm}),\quad t>0.
\end{array}\right.
\end{eqnarray}
Here
$v_{\pm}>0, u_{\pm}, \theta_{\pm}>0, \omega_{\pm}$ are given constants and we assume that $(v_0,u_0, \theta_0,\omega_{0})(\pm\infty)=(v_{\pm},u_{\pm}, \theta_{\pm},\\ \omega_{\pm})$ as compatibility conditions.

Throughout  this paper,  we assume that the pressure $p(v,\theta)$ and the internal energy $e(t,x)$ are given by
\begin{equation}\label{1.3}
p(v,\theta)=\frac{R\theta}{v}=Bv^{-\gamma}\exp\left(\frac{\gamma-1}{R}s\right),\quad e= \frac{R}{\gamma-1}\theta,
\end{equation}
where $s$ is the entropy of the fluid,  and $\gamma>1, B$ and $R$ are positive constants.

The mathematical theory of the compressible micropolar fluid model  has been studied extensively in the last several decades. For the non-isentropic case,  Mujakovi\'{c} first analyzed  the one-dimensional model and obtained a series of results concerning the local-in-time existence and uniqueness, the global existence  and regularity of solutions  to an initial-boundary value problem with homogeneous \cite{N. Mujakovi-1998-1,N. Mujakovi-1998-2,N. Mujakovi-2001} and non-homogeneous  \cite{N. Mujakovi-2007,N. Mujakovi-2008,N. Mujakovi-2009} boundary conditions. Besides, she also studied the large time behavior of the solutions  to initial-boundary value problem  \cite{N. Mujakovi-2005} and the Cauchy problem \cite{N. Mujakovi-2010} of the one-dimensional model. The 1-D compressible viscous micropolar fluid model was also studied by  many other authors, such as Qin et al.  \cite{Qin-2012} proved the stabilization and the regularity of
 solutions  with weighted small initial data,  Duan \cite{Duan-2018-1,Duan-2018} investigated the global existence of strong solutions, and the authors in \cite{Q.Q. Liu-2017,Duan-2017-2,Cui-2017} showed  the  nonlinear stability of some  basic waves (such as rarefaction waves, viscous  contact wave and  viscous  shock wave etc.).  For the 3-D compressible micropolar fluid model, Dra\v{z}i\'{c} and  N. Mujakovi\'{c} considered the local and  global existence, uniqueness, large time behavior and regularity of spherical symmetry solutions, see \cite{I.N-2012,I.N-2015,I.N-2014,I.N-20151,I.N-2016,I.N-2017} and the references therein.

For the isentropic case, Chen \cite{M.T. Chen-2011} proved the global existence of strong solutions to the one-dimensional  model with initial vacuum. Later, Chen and his his collaborators further obtained the global existence of weak solutions \cite{ M.T. Chen-2015} and various  blowup criteria  of strong solutions \cite{M.T. Chen-2012,M.T. Chen-2013} for the 3-D compressible  micropolar fluid system. Liu and Zhang  \cite{Q.Q. Liu-2016} established the optimal time decay of strong solutions to  the 3-D compressible micropolar fluid model. We also mention that there have been many results on the incompressible micropolar
fluid system, see \cite{Q.L. Chen-2012,B. Nowakowski-2013,Dong-2017} and the references therein.

However, few results have been obtained for the compressible  non-viscous  micropolar
fluid model as far. Especially, there no results available now for the nonlinear stability of composite waves for the one-dimensional  non-viscous  micropolar fluid model. This paper is devoted to this problem, and we are concerned with the time-asymptotic nonlinear  stability of the  combination of viscous contact wave with two  rarefaction waves for the Cauchy problem (\ref{1.1})-(\ref{1.2}) when  the far-field states  of the initial data are different, i.e., $(v_+,u_+,\theta_{+},\omega_{+})\neq(v_-,u_-,\theta_{-},\omega_{-})$.

Motivated by the close relationship between the  compressible micropolar fluid model and the compressible Navier-Stokes system as mentioned before, it is expect that the  large-time asymptotic profiles of solutions of the Cauchy problem (\ref{1.1})-(\ref{1.2}) will be the same as those of the compressible Navier-Stokes system in the case of $\omega_{+}=\omega_{-}=0$. It is well-known that the large-time behavior of solutions to the Cauchy problem of the compressible Navier-Stokes system can be described by  the Riemann solution  of the corresponding  Euler system:
\begin{equation}\label{1.4}
\left\{\begin{array}{l}
    v_t-u_x=0, \\[2mm]
    u_t+p(v,\theta)_x=0, \\[2mm]
     \displaystyle\frac{R}{\gamma-1}\theta_t+p(v,\theta)u_x=0, \,\, x\in\mathbb{R},\,\, t>0
  \end{array}\right.
\end{equation}
with the  Riemann  initial data
\begin{equation}\label{1.5}
(v,u,\theta)(0,x)=\left\{\begin{array}{c}
                             (v_-,u_-,\theta_-),\quad   x<0, \\[2mm]
                             (v_+,u_+,\theta_+),\quad   x>0.
                           \end{array}\right.
\end{equation}
Then it is well known \cite{J. Smoller-1994} that  the Euler system (\ref{1.4}) is a strict hyperbolic system of conservation laws with three distinct eigenvalues:
\[\lambda_1(v,\theta)=-\sqrt{\frac{\gamma p}{v}},\quad\lambda_2=0,\quad \lambda_3(v,\theta)=\sqrt{\frac{\gamma p}{v}},\]
and the Riemann solutions of the problem  (\ref{1.4})-(\ref{1.5}) have the wave  patterns: the shock wave, rarefaction wave and contact discontinuity, and  certain  linear combinations  of these three basic waves.  Therefore,  we hope that the solutions of the Cauchy problem (\ref{1.1})-(\ref{1.2}) will tend to some basic waves of the Euler system in the case of $w(t,x)=0$.

In the following two subsections,  we shall construct the viscous contact wave and the combination of viscous contact wave with two  rarefaction waves for the Cauchy problem (\ref{1.1})-(\ref{1.2}), respectively.

\subsection{Viscous contact wave}
It is known that the contact discontinuity solution of the Riemann problem (\ref{1.4})-(\ref{1.5}) takes the form \cite{J. Smoller-1994}
\begin{equation}\label{1.6}
(v^{cd},u^{cd},\theta^{cd})(t,x)=\left\{\begin{array}{c}
                                          (v_-,u_-,\theta_-),\quad   x<0,\ t>0, \\[2mm]
                                          (v_+,u_+,\theta_+),\quad   x>0,\ t<0,
                                        \end{array}\right.
\end{equation}
provided that
\begin{equation}\label{1.7}
u_-=u_+,\qquad p_-\triangleq\frac{R\theta_-}{v_-}=\frac{R\theta_+}{v_+}\triangleq p_+.
\end{equation}

The viscous contact wave $(V^c,U^c,\Theta^c,W^c)(t,x)$ with $W^c(t,x)=0$ corresponding to  the contact discontinuity  $(v^{cd},u^{cd},\theta^{cd},w^{cd})(t,x)$ with $(v^{cd},u^{cd},\theta^{cd})$ defined in (\ref{1.6}) and $w^{cd}(t,x)=0$ for the compressible micropolar fluid model  (\ref{1.1}) becomes smooth and behaviors as a diffusion waves due to the  effect of heat conductivity. As \cite{Q.Q. Liu-2017},  the  viscous contact wave $(V^c, U^c, \Theta^c)(t,x)$ can be defined as follows.

Since the contact discontinuity wave is expected to be almost constant, we set
\begin{equation}\label{1.8}
  P^c=\frac{R\Theta^c}{V^c}= p_+=p_-, \quad |U^c|\ll1,
\end{equation}
then the leading part of the energy equation  $(\ref{1.1})_3$ is
\begin{equation}\label{1.9}
\frac{R}{\gamma-1}\Theta^c_t+p_+U^c_x=\kappa\left(\frac{\Theta^c_x}{V^c}\right)_x.
\end{equation}
Using the equations (\ref{1.8}), $V^c_t=U^c_x$ and (\ref{1.9}),  we get a nonlinear diffusion equation
\begin{equation}\label{1.10}
\left\{\begin{array}{l}
         \Theta^c_t=a\displaystyle\kappa\left(\frac{\Theta^c_x}{\Theta^c}\right)_x,\quad a=\displaystyle\frac{p_+(\gamma-1)}{\gamma R^2}, \\[3mm]
      \Theta^c(\pm\infty,t)=\theta_{\pm}.
       \end{array}\right.
\end{equation}
Due to \cite{C.J. van Duyn-L.A. Peletier-1996/1997}, (\ref{1.10})  has a unique self-similar  solution $\Theta^c(t,x)=\Theta^c(\xi)$ with $\xi=\frac{x}{\sqrt{1+t}}$, which is a monotone function, increasing if $\theta_+>\theta_-$ and decreasing if $\theta_+<\theta_-$. On the other hand, there exists some positive constant $\tilde{\delta}$, such that for $\delta=|\theta_+-\theta_-|\leq\tilde{\delta}$, $\Theta(t,x)$ satisfies
\begin{equation}\label{1.11}
(1+t)^{\frac{k}{2}}\left|\partial_x^k\Theta\right|+|\Theta-\theta_\pm|\leq C\delta e^{-\frac{c_0x^2}{1+t}}, \quad k\geq1,\quad as\,|x|\rightarrow\infty,
\end{equation}
where $c_0$ and $C$ are two positive constants depending only on $\theta_-$ and $\tilde{\delta}$.

Once $\Theta^c(t,x)$ is determined, the viscous contact wave $(V^c, U^c, \Theta^c)(t,x)$ is defined by
\begin{equation}\label{1.12}
  V^c=\frac{R{\Theta^c}}{p_+},\quad U^c=u_-+\frac{\kappa(\gamma-1)}{\gamma R}\frac{\Theta^c_x}{\Theta^c},\quad {\Theta^c}=\Theta^c(t,x),
\end{equation}
then it is easy to check  that  the viscous contact wave $(V^c, U^c, \Theta^c)(t, x)$ satisfies
\begin{eqnarray}\label{1.13}
  \left\{\aligned
          &V^c_t-U^c_x=0, \\[2mm]
          &U^c_t+\displaystyle\left(\frac{R\Theta^c}{V^c}\right)_x=R_1, \\[2mm]
          &\displaystyle\frac{R}{\gamma-1}{\Theta_t^c}+p(V^c, \Theta^c)U^c_x=\displaystyle\kappa\left(\frac{ \Theta^c_x}{V^c}\right)_x
  \endaligned\right.
\end{eqnarray}
with
\begin{equation}\label{1.14}
\aligned
    R_1 &= U_t=\displaystyle O(\delta)(1+t)^{-\frac{3}{2}}e^{-\frac{c_0x^2}{(1+t)}}. \\
    \endaligned
\end{equation}

Our first result is stated as follows.
\begin{Theorem}[Stability of viscous contact wave]
For any given left end  state  $(v_-, u_-, \theta_-)$, suppose that the right end state $(v_+, u_+, \theta_+)$ satisfies (\ref{1.7}).  Let $(V^c, U^c, \Theta^c)(t,x)$ be the viscous contact wave defined in (\ref{1.12}) with strength $\delta=|\theta_+-\theta_-|$.
Then there exist two suitably small positive constants $\epsilon_1$ and $\delta_1\leq(\min\{\tilde{\delta},1\})$ such that if $0<\delta\leq\delta_1$ and
\[
\left\|\left(v_0(\cdot)-V^c(0,\cdot), u_0(\cdot)-U^c(0,\cdot), \theta_0(\cdot)-\Theta^c(0,\cdot), w_0(\cdot)\right) \right\|_{H^2(\mathbb{R})}\leq\varepsilon_1,
\]
then the Cauchy problem (\ref{1.1})-(\ref{1.2}) admits a unique global smooth solution $(v,u,\theta,w)(t,x)$ satisfying
\begin{equation}\label{1.15}
\left\{\aligned
&(v-V^c,u-U^c, \theta-\Theta^c, w)(t,x)\in C(0, +\infty; H^2(\mathbb{R})),\\
&(v-V^c,u-U^c)_x(t,x)\in L^2(0, +\infty; H^{1}(\mathbb{R})), \\
&(\theta-\Theta^c)_x(t,x)\in L^2(0, +\infty; H^{2}(\mathbb{R})),\,\, w(t,x)\in L^2(0, +\infty; H^{3}(\mathbb{R})).
\endaligned\right.
\end{equation}
Moreover, the following large-time behaviors  hold:
\begin{equation}\label{1.16}
\lim_{t\rightarrow+\infty}\sup_{x\in\mathbb{R}}\left\{\left|\left(v-V^c,u-U^c, \theta-\Theta^c, w\right)(t,x)\right|\right\}=0.
\end{equation}
\end{Theorem}

\subsection{Composite waves: combination  of viscous contact wave with rarefaction waves}
When the relation (\ref{1.7}) fails,  the basic theory of hyperbolic systems of conservation laws \cite{J. Smoller-1994} implies  that for any given constant state $(v_-, u_-,\theta_-)$ with $v_->0$, $u_-\in \mathbb{R}$ and $\theta_->0$,  there exists a neighborhood $\Omega(v_-, u_-,\theta_-)$ of $(v_-, u_-,\theta_-)$ such that for any $(v_+, u_+,\theta_+)\in \Omega(v_-, u_-,\theta_-)$, the Riemannn problem (\ref{1.4})-(\ref{1.5}) has a unique solution. In this paper, we only consider the stability of the combination of the viscous contact wave with rarefaction waves. Consequently, we suppose that
\begin{equation}\label{1.17}
(v_+, u_+,\theta_+)\in R_-C_cR_+(v_-, u_-,\theta_-)\subset\Omega(v_-, u_-,\theta_-), \quad |\theta_+-\theta_-|\leq\delta_2,
\end{equation}
where $\delta_2(\leq\tilde{\delta})$ is a positive constant, $R_-,R_+, C_c$ denote  the 1-rarefaction wave curve, 3-rarefaction wave curve, and  the contact wave curve respectively, and
\[
\begin{array}{l} R_-C_cR_+(v_-, u_-,\theta_-)\triangleq\Bigg\{(v, u,\theta)\in\Omega(v_-, u_-,\theta_-)\Bigg|
s\neq s_-,
\\[2mm]
\displaystyle u\geq u_--\int_{v_-}^{e^{\frac{\gamma-1}{R\gamma}(s_--s)}v}\lambda_-(\eta, s_-)\,d\eta,\quad
   u\geq u_--\int_{e^{\frac{\gamma-1}{R\gamma}(s-s_-)}v_-}^{v}\lambda_+(\eta, s)\,d\eta\Bigg\}\end{array}
\]

 It is known from \cite{J. Smoller-1994} that if $\delta_2$ in (\ref{1.17}) is suitably small,  then there exist a positive constant $C=C(\theta_-,\delta_1)$ and a pair of points $(v_-^m, u^m,\theta_-^m)$ and $(v_+^m, u^m,\theta_+^m)$ in $\Omega(v_-, u_-,\theta_-)$ such that
\begin{equation}\label{1.18}
\frac{R\theta_-^m}{v_-^m}=\frac{R\theta_+^m}{v_+^m}\triangleq p^m, \quad |v_{\pm}^m-v_\pm|+|u^m-u_{\pm}|+|\theta_{\pm}^m-\theta_\pm|\leq C|\theta_--\theta_+|.
\end{equation}
Moreover, the states $(v_-^m, u^m,\theta_-^m)$ and $(v_+^m, u^m,\theta_+^m)$ belong to the  1-rarefaction wave curve $R_-(v_-, u_-,\theta_-)$ and  3-rarefaction wave curve $R_+(v_+, u_+,\theta_+)$  respectively, where
\begin{equation}\label{1.19}
R_\pm(v_\pm, u_\pm,\theta_\pm)=\left\{s=s_\pm, u=u_\pm-\int_{v_\pm}^v\lambda_\pm(\eta,s_\pm)\,d\eta,\,\, v>v_\pm\right\}
\end{equation}
with \[ \aligned
&s=\frac{R}{\gamma-1}\ln\frac{R\theta}{A}+R\ln v, \quad s_\pm=\frac{R}{\gamma-1}\ln\frac{R\theta_\pm}{A}+R\ln v_\pm,\\[2mm]
&\lambda_\pm(v,s)=\pm\sqrt{A\gamma v^{-\gamma-1}e^{(\gamma-1)s/R}}.
 \endaligned\]
The contact discontinuity wave curve $C_c$  is defined by
\begin{equation}\label{1.20}
{C_c}(v_-^m, u^m,\theta_-^m)=\left\{(v,u,\theta)(t,x)|u=u^m,\,
p=p^m,\,v\not\equiv v_-^m\right\}.
\end{equation}

The 1-rarefaction wave $(v_-^r, u^r_-,\theta_-^r)(t,x)$ (respectively the  3-rarefaction wave $(v_+^r, u^r_+,\theta_+^r)(t,x)$) connecting $(v_-, u_-,\theta_-)$ and $(v_-^m, u^m,\theta_-^m)$ (respectively $(v_+^m, u^m,\theta_+^m)$ and $(v_+, u_+,\theta_+)$) is a weak solution of the  Euler system (\ref{1.4}) with the Riemann initial data:

\begin{equation}\label{1.21}
(v_\pm^r, u^r_\pm,\theta_\pm^r)(0,x)=\left\{\begin{array}{c}
                             (v_\pm^m,u^m,\theta_\pm^m),\quad   \pm x<0, \\[2mm]
                             (v_\pm,u_\pm,\theta_\pm),\quad   \pm x>0.
                           \end{array}\right.
\end{equation}
The contact discontinuity  wave $(v^{cd},u^{cd},\theta^{cd})(t,x)$ connecting $(v_-^m, u^m,\theta_-^m)$ and $(v_+^m, u^m,\theta_+^m)$  is defined by
\begin{equation}\label{1.22}
(v^{cd},u^{cd},\theta^{cd})(t,x)=\left\{\begin{array}{c}
                                          (v_-^m, u^m,\theta_-^m),\quad   x<0,\ t>0, \\[2mm]
                                          (v_+^m, u^m,\theta_+^m),\quad   x>0,\ t<0.
                                        \end{array}\right.
\end{equation}

To study the stability problem, we need to   construct the smooth approximations of the rarefaction waves $(v_\pm^r, u^r_\pm,\theta_\pm^r)(t,x)$. As \cite{A. Matsumura-K. Nishihara-1986}, we define the smooth approximate rarefaction waves $(V_\pm^r, U^r_\pm,\Theta_\pm^r)(t,x)$ of $(v_\pm^r, u^r_\pm,\theta_\pm^r)(t,x)$ as follows:
\begin{eqnarray}\label{1.23}
\left\{\begin{array}{ll}
        \lambda_{\pm}(V_\pm^r,s_\pm)=\hat{w}_{\pm}(t+1,x),\\[2mm]
        U^{r}_\pm(t,x)=\displaystyle u_\pm-\int_{v_\pm}^{V^{r}(t,x)}\lambda_\pm(\eta, s_\pm)\,d\eta,\\[4mm]
        \Theta_\pm^r=\theta_\pm(v_\pm)^{\gamma-1}(V_\pm^r)^{1-\gamma},
 \end{array}\right.
\end{eqnarray}
where $\hat{w}_-(t,x)$ (respectively $\hat{w}_+(t,x)$) is the solution of the Cauchy  problem of the Burger equation:
\begin{eqnarray}\label{1.24}
\left\{\begin{array}{ll}
         \hat{w}_{t}+ \hat{w}\hat{w}_{x}=0,\,\, x\in\mathbb{R}, t>0,\\[2mm]
 \hat{w}(0,x)=\displaystyle\frac{\hat{w}_r+\hat{w}_l}{2}+\displaystyle\frac{\hat{w}_r-\hat{w}_l}{2} \tanh x
 \end{array}\right.
\end{eqnarray}
with $\hat{w}_l=\lambda_-(v_-,s_-)$ and $\hat{w}_r=\lambda_-(v_-^m,s_-)$ (respectively $\hat{w}_l=\lambda_+(v_+^m,s_+)$ and $\hat{w}_r=\lambda_+(v_+,s_+)$).

Let $(V^{c},  U^{c},  \Theta^{c})(t,x)$ be  the viscous contact wave defined in (\ref{1.11}) with $(v_\pm, u_\pm,\theta_\pm)$ replaced by $(v^m_\pm, u^m_\pm,\theta_\pm^m)(t,x)$ respectively. Set
\begin{eqnarray}\label{1.25}
\left(\begin{array}{ll}
        V\\[1mm]
         U\\[1mm]
     \Theta\\[1mm]
 \end{array}\right)(t,x)=\left(\begin{array}{c}
        V^{r}_-+V^c+V^{r}_+\\[1mm]
        U^{r}_-+U^c+U^{r}_+\\[1mm]
        \Theta^{r}_-+\Theta^c+\Theta^{r}_+\\[1mm]
 \end{array}\right)(t,x)-\left(\begin{array}{c}
        v_{-}^m+v_+^m\\[1mm]
       2u^m\\[1mm]
        \theta_{-}^m+\theta_{+}^m\\[1mm]
 \end{array}\right),
 \end{eqnarray}
then our second main  result  is the following:

\begin{Theorem}[Stability of composite waves]
Suppose that the  end  states $(v_{\pm}, u_{\pm}, \theta_{\pm})$ satisfy (\ref{1.17}) for some  small constant $\delta_2>0$. Let $(V, U, \Theta)(t,x)$ be the superposition  of the viscous contact wave and  the approximate rarefaction waves defined in (\ref{1.25})  with strength $\delta=|\theta_+-\theta_-|$. Then there exist two small positive constants $\varepsilon_2$ and  $\delta_3(\leq\min\{\tilde{\delta},\delta_2,1\})$  such that if   $0<\delta\leq\delta_3$ and the initial data $(v_0, u_0, \theta_0, w_0)(x)$ satisfies
\[
\left\|\left(v_0(\cdot)-V(0,\cdot), u_0(\cdot)-U(0,\cdot), \theta_0(\cdot)-\Theta(0,\cdot), w_0(\cdot)\right) \right\|_{H^2(\mathbb{R})}\leq\varepsilon_2,
\]
then the Cauchy problem (\ref{1.3})-(\ref{1.4}) admits a unique global smooth solution $(v,u,\theta,w)(t,x)$ satisfying
\begin{equation}\label{1.26}
\left\{\aligned
&(v-V,u-U, \theta-\Theta, w)(t,x)\in C(0, +\infty; H^2(\mathbb{R})),\\
&(v-V,u-U)_x(t,x)\in L^2(0, +\infty; H^{1}(\mathbb{R})), \\
&(\theta-\Theta)_x(t,x)\in L^2(0, +\infty; H^{2}(\mathbb{R})),\,\, w(t,x)\in L^2(0, +\infty; H^{3}(\mathbb{R})),
\endaligned\right.
\end{equation}
and  the following large-time behaviors:
\begin{equation}\label{1.27}
\lim_{t\rightarrow+\infty}\sup_{x\in\mathbb{R}}\left(\begin{array}{c}
       \left|\left(v-v_-^r-V^{c}-v_+^r+v_-^m+v_+^m\right)(t,x)\right|\\[2mm]
       \left|\left(u-u_-^r-u_+^r+u^m\right)(t,x)\right|\\[2mm]
        \left|\left(\theta-\theta_-^r-\Theta^{c}-\theta_+^r+\theta_-^m+\theta_+^m\right)(t,x)\right|\end{array}\right)=0,
\end{equation}
\end{Theorem}

Now we make some comments on the analysis of this paper.  Since there is no viscosity in the momentum equation $(\ref{1.1})_2$, the system $(\ref{1.1})$ is less dissipative than the viscous ones considered in the literatures before.  This is the first main  difficulty of this paper. Thus how to derive some suitable a priori estimates for the Cauchy problem (\ref{1.1})-(\ref{1.2}) becomes focus. Especially, how to control the term $\int_0^t\|(\phi_x,\psi_x)(\tau)\|^2_{L^2}d\tau$? If we perform the  energy estimates directly as \cite{Q.Q. Liu-2017,Qin-2012,Duan-2017-2,Cui-2017}, then we can only control $\int_0^t\|\phi_x(\tau)\|^2_{L^2}d\tau$ by $\int_0^t\|\psi_x(\tau)\|^2_{L^2}d\tau$. However, there is no viscosity in our case, one can not get the a priori estimate of $\int_0^t\|\psi_x(\tau)\|^2_{L^2}d\tau$.   Motivated by \cite{Fan-2015}, we derive the a priori estimates directly in the $H^2$ framework and use $\int_0^t\|(\zeta_x,w_x)(\tau)\|^2_{H^1}d\tau$ to control  $\int_0^t\|\phi_x(\tau)\|^2_{L^2}d\tau$ (see \ref{3.49}) for details).  Another difficulty is to control the possible growth of solutions caused by the interaction of waves from different families in the stability of combination of viscous contact waves with  rarefaction waves.  With the help of  a domain  decomposition technique developed in \cite{Huang-2010} and the properties of the approximate rarefaction waves and  viscous contact wave, we can overcome this difficulty and finally get the desired a priori estimates. To the best of our knowledge, this is first work on the stability of composite waves for the compressible  micropolar fluid model.

Before concluding this section, we remark that there are a lot of results on the nonlinear stability of basic waves  for the compressible Navier-Stokes equations. We refer to \cite{A. Matsumura-K. Nishihara-1985,Matsumura-Mei-1999} and the references therein  for the nonlinear  stability of viscous shock waves, \cite{Liu-1988,A. Matsumura-K. Nishihara-1986,Nishihara-2004}  and the references therein for the nonlinear stability of rarefaction waves, \cite{F. M-2005,F. M-2008,Ma-2016}  and the references therein for the nonlinear stability of contact discontinuity, and  \cite{Huang-2016,Huang-Matsumura-2009,Huang-2010,B.-K.-Huang-2017,Fan-2015,Wan-2016} and the references therein for the nonlinear stability of composite waves.

The rest of this paper is organized as follows.  In Section 2, we list some basic properties of the viscous contact wave and rarefaction waves  for later use. An important lemma concerning  the heat kernel and a domain decomposition technique were also presented in this section.  Sections 3 is devoted to the proof of main theorems of this paper. In the final Sections 4, we give the proof of Lemma 3.3 which play an important role in the energy estimates.

{\bf Notations:} Throughout this paper, $C$  stands for some generic positive
constant which may vary in different  estimates. If the dependence need to be explicitly pointed out,
the natation  $C(\cdot,\cdots,\cdot)$ or $C_i(\cdot,\cdots,\cdot)(i\in
{\mathbb{N}})$ is used. For function spaces, $L^p(\mathbb{R})$($1\leq p\leq+\infty$) denotes  the standard  Lebesgue space with
the norm
$\|f\|_{L^p(\mathbb{R})}=\left(\int_\mathbb{R}\left|f(x)\right|^pdx\right)^{\frac{1}{p}}$, and
 $W^{k,p}(\mathbb{R})$ is the usual $k$-th order Sobolev space with its norm
$\|f\|_{W^{k,p}}=\left(\sum_{i=0}^{k}\|\partial_x^if\|^2_{L^p}\right)^{\frac{1}{2}}$. When $p=2$, we simply denote the  the space $W^{k,p}(\mathbb{R})$ by $H^k(\mathbb{R})$, and the norms $\|\cdot\|_{H^k(\mathbb{R})}$ and  $\|\cdot\|_{L^2}$ by $\|\cdot\|_k$ and $\|\cdot\|$, respectively.

\section{Preliminaries}
\setcounter{equation}{0}

First of all, we give the following lemma on the heat kernel which will play an important role in the analysis of this paper, whose proof can be found in \cite{Huang-2010}.

For $\alpha>0$, we define
\begin{equation}\label{2.1}
h(t,x)=(1+t)^{-\frac{1}{2}}\exp\left\{-\frac{\alpha x^2}{1+t}\right\},\quad g(t,x)=\int_{-\infty}^x h(t,y)dy.
\end{equation}
Then it is easy to check that
\begin{equation}\label{2.2}
4\alpha g_t=h_x,\qquad \|g(t)\|_{L^\infty}=\sqrt{\pi}\alpha^{-1/2},
\end{equation}
and  we have
\begin{Lemma}[\cite{Huang-2010}]
For any $0<T\leq\infty$, suppose that the function $F(t,x)$ satisfies
\[F\in L^\infty(0,T; L^2(\mathbb{R})),\quad F_x\in L^2(0,T; L^2(\mathbb{R})),\quad F_t\in L^2(0,T; H^{-1}(\mathbb{R})).\]
Then the following estimate holds:
\begin{equation}\label{2.3}
\int_0^T\int_{\mathbb{R}}F^2h^2dxdt\leq4\pi\|F(0)\|^2+4\pi(\gamma-1)\alpha^{-1}\int_0^T\|F_x(\tau)\|^2d\tau+\frac{8\alpha}{\gamma-1}\int_0^T\langle F_t, Fg^2\rangle\,d\tau,
\end{equation}
where $\langle\cdot,\cdot\rangle$ denotes the inner product on $H^{-1}\times H^1$.
\end{Lemma}

The solution $\hat{w}(t,x)$ of the Cauchy problem (\ref{1.24}) has the following properties.
\begin{Lemma}
For given $\hat{w}_l\in\mathbb{R}$ and $\bar{w}>0$, let $\hat{w}_r\in\{\hat{w}|0<\tilde{w}\triangleq \hat{w}-\hat{w}_l<\bar{w}\}$. Then the problem $(\ref{1.24})$ has a unique global  smooth solution satisfying the following
\begin{itemize}
\item[(i)]$\hat{w}_l<\hat{w}(t,x)<\hat{w}_r$, $\hat{w}_x>0$, $x\in\mathbb{R}, t>0.$
\item[(ii)] For any $p\in[1,+\infty]$, there  exists some positive constant $C=C(p,\hat{w}_l,\bar{w})$ such that for $\tilde{w}\geq0$ and $t\geq0$,
\[\|\hat{w}_x(t)\|_{L^p}\leq C\min\{\tilde{w}, \tilde{w}^{\frac{1}{p}}(1+t)^{-1+\frac{1}{p}}\},\quad
\|\hat{w}_{xx}(t)\|_{L^p}\leq C\min\{\tilde{w}, (1+t)^{-1}\}.\]
\item[(iii)] If $\hat{w}_l>0$, for any $(t,x)\in[0,+\infty)\times(-\infty,0]$,
\[\displaystyle|\hat{w}(t,x)-\hat{w}_l|\leq \tilde{w}e^{-2(|x|+w_lt)},\quad\displaystyle|\hat{w}_x(t,x)|\leq 2\tilde{w}e^{-2(|x|+\hat{w}_lt)}.
\]
\item[(iv)]If $\hat{w}_r<0$, for any $(t,x)\in[0,+\infty)\times[0,+\infty)$,
\[\displaystyle|\hat{w}(t,x)-\hat{w}_r|\leq \tilde{w}e^{-2(|x|+|\hat{w}_r|t)},\quad\displaystyle|\hat{w}_x(t,x)|\leq 2\tilde{w}e^{-2(|x|+|\hat{w}_r|t)}.
\]
\item[(v)] Let $w^{r}(\frac{x}{t})$ be the Riemann solution of the scalar equation $(\ref{1.24})_1$ with the Riemann initial data
\[\hat{w}(0,x)=\left\{\begin{array}{c}
       \hat{w}_l, \quad x<0,\\[2mm]
        \hat{w}_r, \quad x>0,\end{array}
  \right.\]
  then we have \[\lim_{t\rightarrow+\infty}\sup_{x\in\mathbb{R}}\left|\hat{w}(t,x)-w^r\left(\frac{x}{t}\right)\right|=0.\]
\end{itemize}
\end{Lemma}

In order to use Lemma 2.2 to study the properties of the smooth rarefaction waves $(V_{\pm}^r, U_{\pm}^r, \Theta_{\pm}^r)$ constructed in (\ref{1.23}) and the viscous contact wave  $(V^{c}, U^{c}, \Theta^{c})(t,x)$, we divided the the domain $\mathbb{R}\times(0,t)$ into three parts, that is $\mathbb{R}\times(0,t)=\Omega_{-}\cup\Omega_{c}\cup\Omega_{+}$ with
\[
\Omega_{\pm}=\{(x,t)|\pm2x>\pm\lambda_\pm(v_\pm^m,s_\pm)t\},\quad
\Omega_{c}=\{(x,t)|\lambda_-(v_-^m,s_-)t\leq2x\leq\lambda_+(v_+^m,s_+)t\}.
\]

\begin{Lemma}
Assume that (\ref{1.18}) holds with $\delta=|\theta_+-\theta_-|\leq\tilde{\delta}$. Then the smooth rarefaction waves $(V_{\pm}^r, U_{\pm}^r, \Theta_{\pm}^r)$ constructed in (\ref{1.23}) and the viscous contact wave  $(V^{c}, U^{c}, \Theta^{c})(t,x)$  satisfy the following
\begin{itemize}
\item[(i)]$ \left(U_{\pm}^r\right)_x\geq0$, $x\in\mathbb{R},t>0.$
\item[(ii)] For any $p\in[1,+\infty]$, there  exists a positive constant $C=C(p,v_-,u_-, \theta_-, \delta_1, m_0)$ such that for $\delta=|\theta_+-\theta_-|$ and $t\geq0$,
\[ \aligned
&\left\|\left((V_{\pm}^r)_x, (U_{\pm}^r)_x, (\Theta_{\pm}^r)_x\right)(t)\right\|_{L^p}\leq C\min\{\delta, \delta^{\frac{1}{p}}(1+t)^{-1+\frac{1}{p}}\},\\
&\left\|\left(\partial_x^kV_{\pm}^r, \partial_x^kU_{\pm}^r, \partial_x^k\Theta_{\pm}^r\right)(t)\right\|_{L^p}\leq C\min\{\delta, (1+t)^{-1}\}, \quad k=2,3.
 \endaligned\]
\item[(iii)] There exists some positive constant $C=C(p,v_-,u_-, \theta_-, \tilde{\delta}, \delta_2)$ such that for $\delta=|\theta_+-\theta_-|$ and
\[c_1=\frac{1}{10}\min\left\{\left|\lambda_-(v_-^m,s_-)\right|, \lambda_+(v_+^m,s_+), c_0\lambda_-^2(v_-^m,s_-),c_0\lambda_+^2(v_+^m,s_+), 1\right\},\]
we have in $\Omega_{c}$ that
\[
\left|\left((V_{\pm}^r)_x, (U_{\pm}^r)_x, (\Theta_{\pm}^r)_x\right)\right|+\left|V_{\pm}^r-v_{\pm}^m\right|+\left|\Theta_{\pm}^r-\theta_{\pm}^m\right|\leq C\delta e^{-c_1(|x|+t)},
\]
and in  $\Omega_{\mp}$,
\[
\displaystyle\left|V^{c}_x\right|+\left|\Theta^{c}_x\right|+\displaystyle\left|V^{c}-v_{\mp}^m\right|+\left|\Theta^{c}-\theta_{\mp}^m\right|+\left|U^{c}_x\right|\leq\displaystyle C\delta e^{-c_1(|x|+t)},
\]
\[
\left|\left((V_{\pm}^r)_x, (U_{\pm}^r)_x, (\Theta_{\pm}^r)_x\right)\right|+\left|V_{\pm}^r-v_{\pm}^m\right|+\left|\Theta_{\pm}^r-\theta_{\pm}^m\right|\leq C\delta e^{-c_1(|x|+t)}.
\]
\item[(iv)]  It holds that \[\lim_{t\rightarrow+\infty}\sup_{x\in\mathbb{R}}\left|(V_{\pm}^r, U_{\pm}^r, \Theta_{\pm}^r)(t,x)-(v_{\pm}^r, u_{\pm}^r, \theta_{\pm}^r)\left(\frac{x}{t}\right)\right|=0.\]
\end{itemize}
\end{Lemma}

\section{Proof of the main results }
\setcounter{equation}{0}
In this section, we shall prove the time-asymptotic stability of solutions of the Cauchy problem (\ref{1.1})-(\ref{1.2}). Notice that when $(V_{\pm}^m, U_{\pm}^m, \Theta_{\pm}^m)=(V_{\pm}, U_{\pm}, \Theta_{\pm})$, Theorem 1.2 will coincide with Theorem 1.1. Since  Theorem 1.1 can be  proved in the same way as that of Theorem 1.2,  we only give here the proof of Theorem 1.2 for brevity.

\subsection{Reformulation of the  problem}
 First, note that the viscous contact wave $(V^c,U^c,\Theta^c)$ satisfy (\ref{1.13}), and the rarefaction waves $(V_{\pm}^r,U_{\pm}^r,\Theta_{\pm}^r)(x,t)$ solve the Euler equations
\begin{equation}\label{3.1}
  \left\{\begin{array}{l}
    (V_{\pm}^r)_t-(U_{\pm}^r)_x=0, \\[2mm]
    (U_{\pm}^r)_t+p(V_{\pm}^r,\Theta_{\pm}^r)_x=0, \\[2mm]
      \displaystyle\frac{R}{\gamma-1}(\Theta_{\pm}^r)_t+p(V_{\pm}^r,\Theta_{\pm}^r)(U_{\pm}^r)_x=0,
  \end{array}\right.
\end{equation}
it is easy to check that the composite wave $(V,U,\Theta)(x,t)$ defined in (\ref{1.25}) satisfy
 \begin{eqnarray}\label{3.2}
\left\{\begin{array}{ll}
         V_{t}-U_{x}=0,\\[2mm]
         U_{t}+P_{x}=-R_{1},\\[2mm]
        \displaystyle \frac{R}{\gamma-1}\Theta_{t}+PU_{x}=\displaystyle\left(\frac{\kappa\Theta_{x}}{V}\right)_{x}-R_{2},
\end{array}\right.
\end{eqnarray}
where
\[\aligned
R_{1}&=-(P-P_{-}-P_{+})_{x}+U_{t}^{c}:=R_{1}^{1}+U_{t}^{c},\quad P=\frac{R\Theta}{V},\,\,  P_{\pm}=\frac{R\Theta^r_{\pm}}{V^r_{\pm}},\\
R_{2}&=(p^{m}-P)U_{x}^{c}+(P_{-}-P)\left(U_{-}^{r}\right)_x+(P_{+}-P)\left(U_{+}^{r}\right)_x
+\kappa\left(\frac{\Theta_{x}}{V}-\frac{\Theta_{x}^{c}}{V^{c}}\right)_{x}\\
&:=R_{2}^{1}+R_{2}^{2}.
\endaligned\]

 Set the perturbation $(\phi,\psi,\zeta)(t,x)$ by
\[
  \phi(t,x)=v(t,x)-V(t,x),\, \psi(t,x)=u(t,x)-U(t,x),\, \zeta(t,x)=\theta(t,x)-\Theta(t,x),
\]
then it follows from (\ref{1.1}) and (\ref{3.2}) that
\begin{equation}\label{3.3}
  \left\{\begin{array}{ll}
           \phi_t-\psi_x=0, \\[2mm]
           \psi_t+ \displaystyle\left(\frac{R\zeta-P\phi}{v}\right)_x=R_1, \\[2mm]
           \displaystyle\frac{R}{\gamma-1}\zeta_t+p(v,\theta)\psi_x+(p-P)U_x
           =\displaystyle\kappa\left(\frac{\theta_x}{v}-\frac{\Theta_x}{V}\right)_x+\frac{\omega_{x}^{2}}{v}+v\omega^{2}+R_2,\\[2mm]
           \omega_{t}=\displaystyle A\left[\left(\frac{\omega_{x}}{v}\right)_{x}-v\omega\right]
         \end{array}\right.
\end{equation}
 with the following initial  and far-field conditions:
\begin{equation}\label{3.4}
  \left\{\begin{array}{l}
    (\phi,\psi,\zeta)(0,x)=(\phi_0,\psi_0,\zeta_0)(x)=(v-V, u-U, \theta-\Theta)(0,x), \\[2mm]
    (\phi,\psi,\zeta)(t,\pm\infty)=0.
  \end{array}\right.
\end{equation}

 We seek  the solutions of  the Cauchy problem (\ref{3.3})-(\ref{3.4}) in the following set of functions:
\begin{equation}\label{3.5}
\begin{split}
X(0,t)&=\left\{(\phi,\psi,\zeta,\omega)|(\phi,\psi,\zeta,\omega)\in C([0,t],H^{2}(\mathbb{R})),\,\,\phi_{x}\in L^{2}([0,t],H^{1}(\mathbb{R})),\right.
\\
&\left.\quad (\psi,\zeta)_{x}\in L^{2}([0,t],H^{2}(\mathbb{R})),\,\, \omega\in L^{2}([0,t],H^{3}(\mathbb{R}))\right\}.
\end{split}
\end{equation}
where $t>0$ is a positive constant.

The local existence of the Cauchy problem (\ref{3.3})-(\ref{3.4}) is standard, which can be obtained by the iteration technique (See \cite{N. Mujakovi-1998-1}).Thus we omit its proof here for brevity.

To show the global existence of solutions to the Cauchy problem (\ref{3.3})-(\ref{3.4}), it suffice to prove  the following a priori estimates.
\begin{Proposition}[A priori estimates] Under the assumptions of Theorem 1.2, suppose that $(\phi,\psi,\zeta, w)\\(t,x)\in X(0,T)$ for some positive constant $T>0$, and satisfies the following a priori assumption:
\begin{equation}\label{3.6}
\sup_{0\leq t\leq T}\|(\phi,\varphi,\omega,\zeta)(t)\|_{2}\leq\varepsilon
\end{equation}
for some small positive constant $\varepsilon$.
Then there exist two small positive constants $\varepsilon_3$ and  $\delta_3(\leq\min\{\tilde{\delta},\delta_2,1\})$ and a constant  $C_0$ which are independent of $T$,  such that if   $0<\varepsilon\leq\varepsilon_3$ and $0<\delta=|\theta_+-\theta_-|\leq\delta_3$,
then it holds that
\begin{equation}\label{3.7}
\begin{split}
&\|(\phi,\psi,\zeta,\omega)(t)\|_{2}^{2}+\int_{0}^{t}\left(\|(\phi_{x},\psi_{x})(\tau)\|_{1}^{2}+\| \zeta_x(\tau)\|_{2}^{2}+
\|\omega(\tau)\|_{3}^{2}\right)d\tau
\\&\leq C_0\left(\|(\phi_{0},\psi_{0},\zeta_{0},\omega_{0})\|_{2}^{2}+\delta^{\frac{1}{8}}\right)
\end{split}
\end{equation}
for all  $t\in[0,T]$.
\end{Proposition}

\subsection{A priori estimates}
This subsection is devoted to proving Proposition 2.1, which follows a series of Lemmas below.
 First of all, we have from the a priori assumption (\ref{3.6}) and the Sobolev  inequality
\begin{equation}\label{3.8}
\|f\|_{L^\infty}\leq\|f\|^{\frac{1}{2}}\|f_x\|^{\frac{1}{2}},\quad \forall f(x)\in H^1(\mathbb{R})
\end{equation} that
 \begin{equation}\label{3.9}
 \displaystyle\|(\phi, \psi, \zeta, \phi_x, \psi_x, \zeta_x, w, w_x)(t)\|_{{L^\infty}(\mathbb{R})}\leq\varepsilon,\quad \forall\, t\in[0,T].
 \end{equation}
 Furthermore, by the smallness of $\varepsilon$, we have
 \begin{equation}\label{3.10}
 0<\frac{1}{2}\min\{v_-,v_+\}\leq v(t,x)=\phi(t,x)+V(t,x)\leq\frac{3}{2}\max\{v_-,v_+\},
 \end{equation}
and
\begin{equation}\label{3.11}
 0<\frac{1}{2}\min\{\theta_-,\theta_+\}\leq \theta(t,x)=\zeta(t,x)+\Theta(t,x)\leq\frac{3}{2}\max\{\theta_-,\theta_+\}.
\end{equation}

For later use, we first  give the time-decay estimates of the reminder terms $(R_1,R_2)$.
\begin{Lemma}Under the assumption Proposition 3.1, we have
\begin{equation}\label{3.12}
\left\{\begin{array}{ll}
 \|R_{1}(t)\|_{L^\infty}\le C\delta(1+t)^{-\frac{3}{2}},\\[2mm]
 \|R_{2}(t)\|_{L^\infty}\le C(1+t)^{-1},\\[2mm]
          \|R_{1}(t)\|_{W^{2,1}}\le C\delta(1+t)^{-1},\\[2mm]
          \|R_{2}(t)\|_{W^{2,1}}\le C\delta^{\frac{1}{8}}(1+t)^{-\frac{7}{8}},\\[2mm]
         \|R_{1}(t)\|_{2}\le C\delta(1+t)^{-\frac{5}{4}}.
\end{array}\right.
\end{equation}
\end{Lemma}
\noindent{\bf Proof.}~~ First, we have by a direct calculation that
\begin{equation}\label{3.13}
\begin{split}
R_{1}^{1}&=\left(\frac{\Theta_{-}^{r}}{V_{-}^{r}}+\frac{\Theta_{+}^{r}}{V_{+}^{r}}+\frac{\Theta^{c}}{V^{c}}-\frac{\Theta}{V}\right)_{x}
\\&=\left(\frac{\Theta_{-}^{r}}{V_{-}^{r}}+\frac{\Theta_{+}^{r}}{V_{+}^{r}}+\frac{\Theta^{c}}{V^{c}}
-\frac{\Theta_{-}^{r}+\Theta_{+}^{r}+\Theta^{c}-\Theta_{-}^{m}-\Theta_{+}^{m}}{V}\right)_{x}
\\&=(\Theta_{-}^{r})_x\left(\frac{1}{V_{-}^{r}}-\frac{1}{V}\right)+(\Theta_{+}^{r})_x\left(\frac{1}{V_{+}^{r}}-\frac{1}{V}\right)+\Theta_{x}^{c}\left(\frac{1}{V^{c}}-\frac{1}{V}\right)
\\&~~~~+(V_{-}^{r})_x\left(\frac{\Theta}{V^{2}}-\frac{\Theta_{-}^{r}}{(V_{-}^{r})^{2}}\right)+(V_{+}^{r})_x\left(\frac{\Theta}{V^{2}}-\frac{\Theta_{+}^{r}}{(V_{+}^{r})^{2}}\right)+
V^{c}\left(\frac{\Theta}{V^{2}}-\frac{\Theta^{c}}{(V^{c})^{2}}\right).
\end{split}
\end{equation}
It follows from  (\ref{1.25}) and Lemma 2.3 that
\begin{equation}\label{3.14}
  \begin{split}
\left|(\Theta_-^r)_x\left(V^{-1}-(V_-^r)^{-1}\right)\right|
    \le  & C|(\Theta_-^r)_x|\left(|V_+^r-v_+^m|+|V^c-v_-^m|\right)\\
    \le  & C|(\Theta_-^r)_x|\left(|V_+^r-v_+^m|+|V^c-v_-^m|\right)|_{\Omega_-}\\
    & +C|(\Theta_-^r)_x|\left(|V_+^r-v_+^m|+|V^c-v_¡ª^m|\right)|_{\Omega_+ \cup \Omega_c}\\
   \le & C\delta\left[\left(|V_+^r-v_+^m|+|V^c-v_+^m|\right)|_{\Omega_-} +|(\Theta_-^r)_x|_{\Omega_+ \cup \Omega_c}\right]\\
   \le & C\delta^2e^{-c_1(|x|+t)}.
   \end{split}
\end{equation}
Similarly,
\begin{equation}\label{3.15}
 \begin{split}
    & \left|(\Theta_+^r)_x\left(V^{-1}-(V_+^r)^{-1}\right)\right|+\left|(V_{\pm}^r)_x\left(\frac{\Theta_{\pm}^r}{(V_{\pm}^r)^2}-\frac{\Theta}{V^2}\right)\right|\le C\delta^2e^{-c_1(|x|+t)}, \\
     & \left|\Theta_x^c\left(V^{-1}-(V^c)^{-1}\right)\right|+\left|V_x^c\left(\frac{\Theta^c}{(V^c)^2}-\frac{\Theta}{V^2}\right)\right|\le C\delta^{2}e^{-c_1(|x|+t)}.
 \end{split}
\end{equation}
Consequently, it follows from (\ref{3.13})-(\ref{3.15}) and (\ref{1.11}) that
\begin{equation}\label{3.16}
|R_{1}|\leq|R_{1}^{1}|+|U_{t}^{C}|\leq C\delta^2 e^{-c_1(|x|+t)}+C\delta(1+t)^{-\frac{3}{2}}e^{-\frac{c_{0}x^{2}}{1+t}}.
\end{equation}
By the same argument as above, we have for  $R_{2}^{1}$ that
\begin{equation}\label{3.17}
\left|R_{2}^{1}\right|\leq C\delta^2 e^{-c_{1}(|x|+t)}.
\end{equation}
For $R_{2}^{2}$, we have
\begin{equation}\label{3.18}
\aligned
R_{2}^{2}&=\left(\frac{(\Theta_{-}^{r})_x}{V}+\frac{(\Theta_{+}^{r})_x}{V}\right)_x+\left(\frac{\Theta_{x}^{C}}{V}-\frac{\Theta_{x}^{c}}{V^{c}}\right)_x\\
&\leq C\left\{|(\Theta_{-}^{r})_{xx}|+|(\Theta_{-}^{r})_{xx}|+|(\Theta_{-}^{r})_{x}||(V_{-}^{r})_{x}|+|(\Theta_{+}^{r})_{x}||(V_{+}^{r})_{x}|
\right.\\
&\quad\left.+|(\Theta_{-}^{r})_{x}|(|(V_{+}^{r})_{x}|+|V^{c}_{x}|+|(\Theta_{+}^{r})_{x}|(|(V_{-}^{r})_{x}|+|V^{c}_{x}|)\right\}\\
&\quad+C\left \{(|\Theta_{xx}^{c}|+|\Theta_{x}^{c}||V_{x}^{c}|)(|V_{3}^{r}-v_{3}^{m}|+|V_{1}^{r}-v_{-}^{m}|)+
|\Theta_{x}^{c}|(|(V_{-}^{r})_x|+|(V_{+}^{r})_x|)\right\}\\
&:=R_{2,1}^{1}+R_{2,2}^2,
\endaligned
\end{equation}
then Lemma 2.3 implies that
\begin{equation}\label{3.19}
\|R_{2,1}^{2}\|_{L^1}\leq C\delta^{\frac{1}{8}}(1+t)^{-\frac{7}{8}},\quad|R_{2,1}^{2}|\leq C (1+t)^{-1}, \quad |R_{2,2}^{2}|\leq C\delta^2 e^{-c_1(|x|+t)}.
\end{equation}
Thus we have  from (\ref{3.16})-(\ref{3.19}) that
\begin{equation}\label{3.191}
\aligned &\|R_{1}(t)\|_{L^\infty}\le C\delta(1+t)^{-\frac{3}{2}},\quad \|R_{2}(t)\|_{L^\infty}\le C(1+t)^{-1},\\
&\|R_{1}(t)\|_{L^1}\le C\delta(1+t)^{-1},\quad\|R_{2}(t)\|_{L^1}\le C\delta^{\frac{1}{8}}(1+t)^{-\frac{7}{8}},\quad\|R_{1}(t)\|_{L^2}\le C\delta(1+t)^{-\frac{5}{4}}.
\endaligned
\end{equation}
Similarly, we can also obtain
\begin{equation}\label{3.192}
\aligned
&\|(R_{1x},R_{1xx})(t)\|_{L^1}\le C\delta(1+t)^{-1},\quad\|(R_{2x},R_{2xx})(t)\|_{L^1}\le C\delta^{\frac{1}{8}}(1+t)^{-\frac{7}{8}},\\
&\|(R_{1x},R_{1xx})(t)\|_{L^2}\le C\delta(1+t)^{-\frac{5}{4}}.
\endaligned
\end{equation}
(\ref{3.12}) thus follows from (\ref{3.191}) and (\ref{3.192}) immediately. This completes the proof of Lemma 3.1.

For the $L^2$ estimates on  $\|(\phi,\psi,\zeta,\omega)(t)\|$, we have
\begin{Lemma}Under the assumption Proposition 3.1, there exists a positive constant $C$ such that
\begin{equation}\label{3.20}
\begin{split}
&\|(\phi,\psi,\zeta,\omega)(t)\|^{2}+\int_{0}^{t}\int_{\mathbb{R}}\left[\left((U_{-}^{r})_x+(U_{+}^{r})_x\right)(\phi^{2}+\zeta^{2})+\zeta_{x}^{2}+\omega_{x}^{2}+\omega^{2}\right]dxd\tau
\\
&\leq C\left(\|(\phi_{0},\psi_{0},\zeta_{0},\omega_{0})\|^{2}+
\delta^{\frac{1}{8}}+\delta^{\frac{1}{8}}\int_{0}^{t}\|(\phi_{x},\psi_{x})(\tau)\|^{2}d\tau\right).
\end{split}
\end{equation}
\end{Lemma}
\noindent{\bf Proof.}~~Multiplying $(\ref{3.3})_{1}$ by $-R\Theta\left(\frac{1}{v}-\frac{1}{V}\right)$, $(\ref{3.3})_{2}$ by $\psi$, $(\ref{3.3})_{3}$ by $\frac{\zeta}{\theta}$, $(\ref{3.3})_{4}$ by $\omega$, and adding the  resultant equations together, we have
\begin{equation}\label{3.21}
\begin{split}
&\left\{\Theta\Phi\left(\frac{v}{V}\right)+\frac{1}{2}\varphi^{2}
+\Theta\Phi\left(\frac{\theta}{\Theta}\right)+\frac{\omega^{2}}{2}\right\}_{t}+\frac{\kappa\zeta_{x}^{2}}{v\theta}+Av\omega^{2}+\frac{A}{v}\omega_{x}^{2}
+H_{1x}+Q_{1}
\\&=R_{1}\psi+R_{2}\frac{\zeta}{\theta}+Q_{2}+Q_{3},
\end{split}
\end{equation}
where
\begin{equation}\label{3.22}
\begin{split}
&H_{1}=(p-P)\varphi-\frac{\kappa\zeta}{\theta}\left(\frac{\zeta_{x}}{v}-\frac{\Theta_{x}\phi}{vV}\right)-\frac{Aww_{x}}{v},\\
&Q_{1}=-R\Theta_{t}\Phi\left(\frac{v}{V}\right)+\frac{PU_{x}}{vV}\phi^{2}+\frac{R\Theta_{t}}{\gamma-1}\Phi\left(\frac{\Theta}{\theta}\right)+\frac{\zeta}{\theta}(p-P)U_{x},\\
&Q_{2}=\kappa\frac{\Theta_{x}}{\theta^{2}v}\zeta\zeta_{x}+\kappa\frac{\zeta_{x}\phi}{\theta vV}\Theta_{x}-\kappa\frac{\theta_{x}\zeta\phi}{\theta^{2}vV}\Theta_{x},\quad Q_{3}=\frac{\zeta}{\theta}\left(\frac{\omega_{x}^{2}}{v}+v\omega^{2}\right),\\
&\Phi(s)=s-1-\ln s.
\end{split}
\end{equation}

Since (\ref{3.2}) implies that
\begin{equation}\label{3.23}
\begin{split}
-R\Theta_{t}&=(\gamma-1)P_{-}(U_{-}^{r})_x+(\gamma-1)P_{+}(U_{-}^{r})_x+(\gamma-1)\left(p^{m}U_{x}^{c}-\left(\frac{\Theta_{x}^{c}}{V^{c}}\right)_{x}\right)\\
&=(\gamma-1)P((U_{-}^{r})_x+(U_{+}^{r})_x)+(\gamma-1)(P_--P)(U_{-}^{r})_x\\
&\quad+(\gamma-1)(P_+-P)(U_{+}^{r})_x-p^{m}U_{x}^{c},
\end{split}
\end{equation}
we have
\begin{equation}\label{3.24}
\begin{split}
Q_{1}&=-R\Theta_{t}\left[\Phi\left(\frac{v}{V}\right)-\frac{1}{\gamma-1}\Phi\left(\frac{\Theta}{\theta}\right)\right]+\left[\frac{P\phi^{2}}{vV}+\frac{\zeta}{\theta}(p-P)\right]U_{x}
\\&=((U_{-}^{r})_x+(U_{+}^{r})_x) Q_{11}+Q_{12},
\end{split}
\end{equation}
where
\begin{equation}\label{3.25}
\begin{split}
Q_{11}&=(\gamma-1)P(V,\Theta)\Phi\left(\frac{v}{V}\right)+\frac{P\phi^2}{vV}-P\Phi\left(\frac{\Theta}{\theta}\right)+\frac{\zeta}{\theta}(p-P)\\
     &= P\left(\Phi\left(\frac{\theta V}{\Theta v}\right)+\gamma\Phi\left(\frac{v}{V}\right)\right)\geq C(\phi^{2}+\zeta^{2}),
\end{split}
\end{equation}
\begin{equation}\label{3.26}
\begin{split}
Q_{12}&=-U_x^c\left[\frac{P(V,\Theta)\phi^2}{vV}-p^m\Phi\left(\frac{v}{V}\right)+\frac{p^m}{\gamma-1}\Phi\left(\frac{\Theta}{\theta}\right)+\frac{\zeta}{\theta}(p-P)\right]
\\&\quad+(\gamma-1)(P-P_-)(U_-^r)_x\left(\Phi\left(\frac{v}{V}\right)-\frac{1}{\gamma-1}\Phi\left(\frac{\Theta}{\theta}\right)\right)
\\&\quad+(\gamma-1)(P-P_+)(U_+^r)_x\left(\Phi\left(\frac{v}{V}\right)-\frac{1}{\gamma-1}\Phi\left(\frac{\Theta}{\theta}\right)\right).
\end{split}
\end{equation}
It follows from (\ref{1.11}), (\ref{3.9}), Lemma 2.3 and the Cauchy inequality that
\begin{equation}\label{3.27}
\aligned
&|Q_{12}|+|Q_{2}|+|Q_{3}|\\
\leq& \eta\zeta_{x}^{2}+C_\eta\Theta_{x}^{2}(\phi^{2}+\zeta^{2})+\varepsilon(w^2+w_x^2)+C\left(|U_{x}^{c}|+\delta((U_{-}^{r})_x+(U_{+}^{r})_x)\right)(\phi^{2}+\zeta^{2})\\
\leq& \eta\zeta_{x}^{2}+\varepsilon(w^2+w_x^2)+C_\eta\delta((U_{-}^{r})_x+(U_{+}^{r})_x)(\phi^{2}+\zeta^{2})
+C\delta(1+t)^{-1}e^{-\frac{c_{1}x^{2}}{1+t}}(\phi^{2}+\zeta^{2}).
\endaligned
\end{equation}
Here and hereafter, $\eta$ is a small positive constant and $C_\eta>0$ is a constant depending on $\eta$,  and in (\ref{3.27}), we have used the fact that
\begin{equation*}\aligned
  |(\Theta_-^r)_x|=&|\theta_-(v_-)^{\gamma-1}(1-\gamma)(V_-^r)^{-\gamma}(V_-^r)_x|\\
  =&\left|\theta_-(v_-)^{\gamma-1}(1-\gamma)(V_-^r)^{-\gamma}\frac{(U_-^r)_x}{\lambda_-(V_-^r,s_-)}\right|\le C(U_-^r)_x.
\endaligned\end{equation*}
Moreover, we derive  from the Sobolev inequality, the Young inequality and Lemma 3.1 that
\begin{equation}\label{3.28}
  \begin{split}
     &\left|\int_0^t\int_{\mathbb R}\left(R_{1}\psi+R_{2}\frac{\zeta}{\theta}\right)dxd\tau\right| \\
     \le&C\int_0^t\|(\psi,\zeta)(\tau)\|^{\frac{1}{2}}\|(\psi_x,\zeta_x)(\tau)\|^{\frac{1}{2}}\|(R_{1},R_{2})(\tau)\|_{L^1}d\tau\\
      \le& \delta^{\frac{1}{8}}\int_0^t\|(\psi_x,\zeta_x)(\tau)\|^2d\tau+C\int_0^t\sup_{0\leq\tau\leq t}\|(\phi,\zeta)(\tau)\|^{\frac{2}{3}} \delta^{\frac{1}{8}}(1+\tau)^{-\frac{7}{6}}d\tau\\
     \le &\delta^{\frac{1}{8}}\int_0^t\|(\psi_x,\zeta_x)(\tau)\|^2d\tau+C\delta^{\frac{1}{8}}.
  \end{split}
\end{equation}

Integrating (\ref{3.21}) in $t$ and $x$ over $[0,t]\times\mathbb{R}$, and using (\ref{3.27}), (\ref{3.28}), the Lemma 3.3 below and the smallness of $\eta,\varepsilon$ and $\delta$, we can get (\ref{3.20}). This completes  the proof of Lemma 3.2.

\begin{Lemma}Under the assumption Proposition 3.1, there exist positive constants $c_2=\frac{c_1}{4}$ and $C$ such that
\begin{equation}\label{3.29}
\begin{split}
&\int_{0}^{t}\int_{\mathbb{R}}(1+\tau)^{-1}e^{-\frac{c_{2}x^{2}}{1+\tau}}(\phi^{2}+\psi^{2}+\zeta^{2})dxd\tau
\\&\leq C+C\int_{0}^{t}\|(\phi_{x},\psi_{x},\zeta_{x}, w_x,w)(\tau)\|^{2}d\tau+C\int_{0}^{t}\int_{\mathbb{R}}\left((U_{-}^{r})_x+(U_{+}^{r})_x\right)(\phi^{2}+\zeta^{2})dxd\tau.
\end{split}
\end{equation}
\end{Lemma}
\noindent{\bf Proof.}~~ The proof of  (\ref{3.29}) is divided into the following two inequalities:
\begin{equation}\label{4.1}
\begin{split}
\int_{0}^{t}\int_{\mathbb{R}}h^{2}[(R\zeta-P\phi)^{2}+\psi^{2}]dxd\tau
&\leq C+C\int_{0}^{t}\left(\|(\phi_{x},\psi_{x},\zeta_{x}, \omega_{x})(\tau)\|^{2}+\|\omega(\tau)\|^2\right)d\tau
\\&\quad +C\delta\int_{0}^{t}\int_{\mathbb{R}}h^{2}(\psi^{2}+\zeta^{2})dxd\tau,
\end{split}
\end{equation}
and for any $\eta>0$,
\begin{equation}\label{4.2}
\begin{split}
\int_{0}^{t}\int_{\mathbb{R}}(R\zeta+(\gamma-1)P\phi)^{2}h^2dxd\tau
&\leq C+C_\eta\int_{0}^{t}\left(\|(\phi_{x},\psi_{x},\zeta_{x}, \omega_{x})(\tau)\|^{2}+\|\omega(\tau)\|^2\right)d\tau
\\&\quad~~~+C(\delta+\eta)\int_{0}^{t}\int_{\mathbb{R}}h^{2}(\psi^{2}+\zeta^{2})dxd\tau
\\&\quad~~~+C\int_{0}^{t}\int_{\mathbb{R}}(\phi^2+\zeta^2)((U_{-}^{r})_x+(U_{+}^{r})_x)dxd\tau.
\end{split}
\end{equation}

Indeed, if we add (\ref{4.2}) to (\ref{4.1}), then taking first $\eta$ and then $\delta$ suitably small, and noticing  that
\[
  [R\zeta+(\gamma-1)P\phi]^2+(R\zeta-P\phi)^2\geq R^2\zeta^2+2(\gamma-1)P^2\phi^2
  \geq c_2\left(\phi^2+\zeta^2\right)
\]
with $c_2$ being a positive constant depending only on $R,\gamma,\theta_+,\theta_-,v_+$ and $v_-$, we can get (\ref{3.29}) immediately.

We first prove (\ref{4.1}). To do so, denote
\[
f(t,x)=\int_{-\infty}^{x}h^{2}(t,y)dy,
\]
then we have
\begin{equation}\label{4.3}
\|f(t)\|_{L^{\infty}}\leq2\alpha^{-\frac{1}{2}}(1+t)^{-\frac{1}{2}},\quad
\|f_{t}(t)\|_{L^{\infty}}\leq4\alpha^{-\frac{1}{2}}(1+t)^{-\frac{3}{2}}
\end{equation}
Multiplying  $(\ref{3.3})_2$ by $(R\zeta-P\phi)vf$, and integrating the resulting equation over $\mathbb{R}$ leads to
\begin{equation}\label{4.4}
\begin{split}
&\frac{1}{2}\int_{\mathbb{R}}(R\zeta-P\phi)^2h^2 dx
\\&=
\int_{\mathbb{R}}\psi_{t}(R\zeta-P\phi)vf dx-\int_{\mathbb{R}}R_{1}(R\zeta-P\phi)vf dx
\\&=
\left(\int_{\mathbb{R}}\psi(R\zeta-P\phi)vf dx\right)_{t}-\int_{\mathbb{R}}\psi(R\zeta-P\phi)_tvf dx-\int_{\mathbb{R}}\psi(R\zeta-p\phi)v_tf dx
\\&\quad -\int_{\mathbb{R}}\psi(R\zeta-P\phi)vf_tdx
-\int_{\mathbb{R}}v^{-1}(R\zeta-P\phi)^2v_xfdx
-\int_{\mathbb{R}}R_{1}(R\zeta-P\phi)vf dx
\\&:=\left(\int_{\mathbb{R}}\psi(R\zeta-P\phi)vf dx\right)_{t}+\sum_{i=1}^{5}L_{i}.
\end{split}
\end{equation}
We now turn to estimate $L_{i}(1\leq i\leq 5)$ term by term. Using equations $(\ref{3.3})_1$ and $(\ref{3.3})_2$, we have
\begin{equation}\label{4.5}
\left(\frac{R}{\gamma-1}\zeta+P\phi\right)_t=-\frac{R\zeta-P\phi}{v}(\psi_x+U_x)+\kappa\left(\frac{V\zeta_x-\phi\Theta_x}{vV}\right)_x+P_t\phi
+v\omega^2+\frac{\omega_{x}^{2}}{v}+R_2,
\end{equation}
thus it holds that
\begin{eqnarray}\label{4.6}
L_{1}&=&-(\gamma-1)\int_{\mathbb{R}}\psi vf\left(\frac{R}{\gamma-1}\zeta+P\phi\right)_tdx+\gamma\int_{\mathbb{R}}(P\phi)_t vf\psi dx\nonumber\\
&=&(\gamma-1)\int_{\mathbb{R}}\psi f(R\zeta-P\phi)(U_{x}+\psi_{x})dx+\kappa(\gamma-1)\int_{\mathbb{R}}\frac{V\zeta_{x}-
\phi\Theta_{x}}{Vv}(\psi vf)_{x}dx\nonumber\\
&&-(\gamma-1)\int_{\mathbb{R}}\psi vf\left(v\omega^2+\frac{\omega_{x}^{2}}{v}+R_{2}\right)dx+\frac{\gamma}{2}\int_\mathbb{R}P vf(\psi^2)_xdx
+\int_{\mathbb{R}}\psi vfP_t\phi dx\nonumber\\
&:=&\sum_{i=1}^{5}L_{1}^{i}.
\end{eqnarray}

 To control  $L_{1}^{i}(1\leq i\leq5)$ and $L_{i}(1\leq i\leq5)$, we derive  from the Cauchy inequality, the Sobolev inequality (\ref{3.8}), the Young inequality, $(\ref{1.11}), (\ref{3.3})_1, (\ref{3.9})$, (\ref{3.12}) and (\ref{4.3}) that
\begin{eqnarray}\label{4.7}
|L_{1}^{1}+L_{2}|&\leq&C\left|\int_{\mathbb{R}}\psi f(R\zeta-P\phi)(U_{x}+\psi_{x})dx\right|\nonumber\\
&\leq& C\|\psi(t)\|_{L^{\infty}}\|f(t)\|_{L^{\infty}}\int_{\mathbb{R}}(|\zeta|+|\phi|)(|U_{x}|+|\psi_{x}|)dx\nonumber\\
&\leq& C\|\psi(t)\|^{\frac{1}{2}}\|\psi_x(t)\|^{\frac{1}{2}}\|f(t)\|_{L^{\infty}}\|(\zeta,\phi)(t)\|(\|((U_{-}^r)_{x},(U_{+}^r)_{x},U_{x}^c)(t)\|+\|\psi_{x}(t)\|)
\nonumber\\&\leq& C\|\psi_{x}(t)\|^{2}+C(1+t)^{-\frac{4}{3}},\\[2mm]
|L_{1}^{2}|&\leq&\int_{\mathbb{R}}(|\zeta_{x}|+|\phi|(1+t)^{-\frac{1}{2}})(|\psi_{x} vf|+|\psi v_{x}f|+|\psi vf_{x}|)dx
\nonumber\\&\leq& C\left(\|\zeta_{x}(t)\|+\|\phi(t)\|(1+t)^{-\frac{1}{2}}\right)
\left(\|(\psi_{x},\phi_x)(t)\|(1+t)^{-\frac{1}{2}}+\|\psi(t)\|(1+t)^{-1}\right)
\nonumber\\&\leq& C\|(\phi_{x},\psi_{x},\zeta_{x})(t)\|^2+C(1+t)^{-\frac{3}{2}},\\[2mm]
|L_{1}^{3}|&\leq& C\|f(t)\|_{L^{\infty}}\int_{\mathbb{R}}\left(\omega^{2}+\omega_{x}^{2}+|R_2|\right)dx
\leq C\|\omega(t)\|_{1}^{2}+C(1+t)^{-\frac{11}{8}},\\[2mm]
|L_{1}^{4}|&\leq& C\int_{\mathbb{R}}f^{\frac{3}{2}}\phi^2\psi^2dx+C\int_{\mathbb{R}}f^{\frac{1}{2}}P_{t}^2dx
\nonumber\\&\leq&\|\psi(t)\|_{L^\infty}^2(1+t)^{-\frac{3}{4}}\|\phi(t)\|^2+C(1+t)^{-\frac{1}{4}}\|((U_{-}^r)_{x},(U_{+}^r)_{x},U_{x}^c)(t)\|^2
\nonumber\\&\leq&C\|\psi(t)\|\|\psi_x(t)\|(1+t)^{-\frac{3}{4}}+C(1+t)^{-\frac{5}{4}}
\nonumber\\&\leq&C\|\psi_x(t)\|^2+C(1+t)^{-\frac{5}{4}},\\[2mm]
L_{1}^{5}&=&-\frac{\gamma P}{2}\int_{\mathbb{R}}\psi^2f(V_x+\phi_x)dx-\frac{\gamma P}{2}\int_{\mathbb{R}}v\psi^2h^2dx-\frac{\gamma}{2}\int_{\mathbb{R}}P_xvf\psi^2dx
\nonumber\\&\leq&-\frac{\gamma P}{2}\int_{\mathbb{R}}v\psi^2h^2dx+C\|f(t)\|_{L^{\infty}}\int_{\mathbb{R}}\psi^2\left(|V_x^c|+\|((V_+^r)_x,(V_-^r)_x)(t)\|_{L^\infty}\right)dx\nonumber\\
&&\quad+C\|f(t)\|_{L^{\infty}}\|\psi(t)\|^{\frac{3}{2}}\|\psi_x(t)\|^{\frac{1}{2}}\|\phi_x(t)\|
\nonumber\\&\leq&-\frac{\gamma P}{2}\int_{\mathbb{R}}v\psi^2h^2dx+C\delta\int_{\mathbb{R}}\psi^2h^2dx+C(1+t)^{-\frac{3}{2}}
\nonumber\\&\leq& -\frac{\gamma p}{4}\int_{\mathbb{R}}vh^{2}\psi^{2}dx+C\|(\phi_x,\psi_x)(t)\|^2+C(1+t)^{-\frac{3}{2}},\\[2mm]
|L_{3}|&\leq& C\|\psi(t)\|(\|\phi(t)\|+\|\zeta(t)\|)\|f_t(t)\|_{L^\infty}
\leq C(1+t)^{-\frac{3}{2}},\\[2mm]
|L_{4}|&\leq& C\left|\int_{\mathbb{R}}(R\zeta-P\phi)^2\left(|(\Theta_{-}^{r})_x|+|(\Theta_{+}^{r})_x|+|\Theta_{x}^{c}|+|\phi_x|\right)fdx\right|
\nonumber\\&\leq& C\|f(t)\|_{L^{\infty}}\int_{\mathbb{R}}(\phi^2+\zeta^2)\left(|\Theta_{x}^{c}|+\left\|((\Theta_{-}^{r})_x,(\Theta_{+}^{r})_x)(t)\right\|_{L^\infty}\right)dx\nonumber\\
&&\quad+C\|f(t)\|_{L^{\infty}}\left(\|\zeta(t)\|^{\frac{3}{2}}\|\zeta_x(t)\|^{\frac{1}{2}}\|\phi_x(t)\|
\|\phi(t)\|^{\frac{3}{2}}\|\phi_x(t)\|^{\frac{3}{2}}\right)
\nonumber\\&\leq& C\delta\int_{\mathbb{R}}(\phi^2+\zeta^2)h^2dx+C(\|\zeta_x(t)\|^2+\|\phi_x(t)\|^2)+C(1+t)^{-2},
\end{eqnarray}
and
\begin{equation}\label{4.14}
|L_{5}|\leq C(1+t)^{-\frac{1}{2}}(\|\zeta(t)\|_{L^{\infty}}+\|\phi(t)\|_{L^{\infty}})\|R_{1}\|_{L^{1}}
\leq C(1+t)^{-\frac{3}{2}}.
\end{equation}

 Integrating $(\ref{4.4})$ over $(0,t)$, and using  the estimates $(\ref{4.6})-(\ref{4.14})$, we have (\ref{4.1}) holds.

Next we prove  (\ref{4.2}). Let $\Gamma=R\zeta+(\gamma-1)P\phi$, then from (\ref{4.5}), we have
\begin{eqnarray}
\langle \Gamma_{t},\Gamma g^{2}\rangle_{H^{-1}\times H^{1}}&=&\int_{\mathbb{R}}(R\zeta+(\gamma-1)P\phi)_{t}\Gamma g^{2}dx
\nonumber\\&=&-\int_\mathbb{R}\frac{\gamma-1}{v}(R\zeta-P\phi)\Gamma\psi_xg^2dx
-\int_\mathbb{R}\frac{\gamma-1}{v}(R\zeta-P\phi)\Gamma U_xg^2dx
\nonumber\\&&\quad-\kappa(\gamma-1)\int_\mathbb{\mathbb{R}}\frac{V\zeta_x-\phi\Theta_x}{vV}(\Gamma g^{2})_xdx+
(\gamma-1)\int_{\mathbb{R}}P_t\phi \Gamma g^2dx
\nonumber\\&&\quad+(\gamma-1)\int_{\mathbb{R}}(R_{2}+v\omega^{2}+\frac{\omega_{x}^{2}}{v})\Gamma g^{2}dx:=\sum_{i=6}^{10}L_i.\nonumber
\end{eqnarray}
Noticing that
$$|U_x|\leq C\left|\left(\frac{\Theta_x}{\Theta}\right)_x\right|\leq C(|\Theta_{xx}|+\Theta_{x}^{2})\leq C\left(\delta h^2+\left|((\Theta_+^r)_{xx},(\Theta_-^r)_{xx})\right|+\left|((\Theta_+^r)_{x},(\Theta_-^r)_{x})\right|^2\right),$$
we derive
\begin{equation}\label{4.15}
\begin{split}
|L_{7}|&\leq C\int_{\mathbb{R}}(\phi^{2}+\zeta^{2})\delta h^{2}dx+C\|(\phi,\zeta)(t)\|\|(\phi_x,\zeta_x)(t)\|\|((\Theta_+^r)_{xx},((\Theta_+^r)_{x})^2)(t)\|_{L^1}\\
&\leq C\|(\phi_x,\zeta_x)(t)\|^2+C\delta\int_{\mathbb{R}}(\phi^{2}+\zeta^{2})h^{2}dx+C(1+t)^{-2}.
\end{split}
\end{equation}
Using ($\ref{1.11}$), Lemma 2.3 and the Cauchy inequality, we get
\begin{equation}\label{4.16}
\begin{split}
|L_{8}|&\leq C\int_{\mathbb{R}}|v^{-1}\zeta_{x}\Gamma_{x}g^{2}|dx+C\int_{R}|v^{-1}\zeta_{x}\Gamma gh|dx
\\&~~~~+C\int_{\mathbb{R}}|v^{-1}\phi\Theta_x\Gamma_xg^2|dx+C\int_{\mathbb{R}}|v^{-1}\phi\Theta_x\Gamma gh|dx
\\&\leq C(\delta+\eta)\int_{\mathbb{R}}(\phi^{2}+\zeta^{2})h^{2}dx+\frac{C}{\eta}(\|\zeta_{x}(t)\|^{2}+\|\phi_{x}(t)\|^{2})+C(1+t)^{-2},
\end{split}
\end{equation}
\begin{eqnarray}\label{4.17}
|L_{9}|&\leq& C\int_{\mathbb{R}}|\phi||\Gamma|\left(|U_{x}^{c}|+(U_{-}^{r})_x+(U_{+}^{r})_x\right)
\nonumber\\
&\leq& C\int_{\mathbb{R}}(\phi^2+|\phi\zeta|)\left(|U_{x}^{c}|+(U_{-}^{r})_x+(U_{+}^{r})_x\right)dx
\nonumber\\
&\leq& C\delta\int_{\mathbb{R}}(\phi^{2}+\zeta^{2})h^{2}dx+C\int_{\mathbb{R}}(\phi^{2}+\zeta^{2})((U_{-}^{r})_x+(U_{+}^{r})_x)dx,
\end{eqnarray}
\begin{equation}\label{4.18}
\begin{split}
|L_{10}|&\leq C\int_{\mathbb{R}}\left(\omega^{2}+\omega_{x}^{2}+|R_{2}|\right)|\Gamma|dx\\
&\leq C\|\omega(t)\|_{1}^{2}+C\|(\zeta,\phi)(t)\|^{\frac{1}{2}}\|(\zeta_x,\phi_x)(t)\|^{\frac{1}{2}}\|R_{2}(t)\|_{L^1}
\\
&\leq C\|\omega(t)\|_{1}^{2}+C\|(\zeta_x,\phi_x)(t)\|^{2}+C(1+t)^{-\frac{7}{6}}.
\end{split}\end{equation}

Finally, we estimate $L_6$. Noticing that $R\zeta-P\phi=H-\gamma P\phi$, we obtain by $(\ref{3.3})_1$, $(\ref{2.2})$ and ($\ref{4.5}$) that
\begin{eqnarray}
\frac{-2L_{6}}{\gamma-1}&=&2\int_\mathbb{R}\frac{1}{v}(\Gamma-\gamma P\phi)\Gamma\psi_xg^2dx
\nonumber\\&=&\int_{\mathbb{R}}\left(2v^{-1}\Gamma^2\phi_tg^2-\gamma Pv^{-1}\Gamma(\phi^{2})_tg^2)\right)dx
\nonumber\\&=&\left(\int_{\mathbb{R}}\frac{2\Gamma-\gamma P\phi}{v}\Gamma g^{2}\phi dx\right)_{t}
-2\int_{\mathbb{R}}\frac{2\Gamma-\gamma P\phi}{v}\Gamma gg_t\phi dx
\nonumber\\&&~~~~+\int_{\mathbb{R}}\frac{2\Gamma-\gamma P\phi}{v^2}\Gamma g^{2}v_t\phi dx-\int_{\mathbb{R}}\frac{4\Gamma-\gamma P\phi}{v}\Gamma_tg^{2}\phi dx+\int_{\mathbb{R}}v^{-1}\gamma P_t\phi^2\Gamma g^2 dx
\nonumber\\&=&\left(\int_{\mathbb{R}}\frac{2\Gamma-\gamma P\phi}{v}\Gamma g^{2}\phi dx\right)_{t}
-\frac{1}{2\alpha}\int_{\mathbb{R}}\frac{2\Gamma-\gamma P\phi}{v}\Gamma gh_x\phi dx
\nonumber\\&&~~~~+\int_{\mathbb{R}}v^{-2}u_xg^{2}\phi\left[\Gamma(2\Gamma-\gamma P\phi)+(\gamma-1)(4\Gamma-\gamma P\phi)(R\zeta-P\phi)\right]dx
\nonumber\\&&~~~~+\kappa(\gamma-1)\int_{\mathbb{R}}\frac{V\zeta_{x}-\phi\Theta_{x}}{vV}\left(\frac{g^{2}\phi(4\Gamma-\gamma P\phi)}{v}\right)_{x}dx
\nonumber\\&&~~~~-(\gamma-1)\int_{\mathbb{R}}\frac{g^{2}\phi}{v}(4\Gamma-\gamma P\phi)(v\omega^{2}+\frac{\omega_{x}^{2}}{v}+R_{2})dx\nonumber\\
&&\quad+\int_{\mathbb{R}}v^{-1}P_t\phi^2g^2[(4-3\gamma)\Gamma+\gamma(\gamma-1)P\phi]dx
\nonumber\\&:=&\left(\int_{\mathbb{R}}\frac{2\Gamma-\gamma P\phi}{v}\Gamma g^{2}\phi dx\right)_{t}+\sum_{i=1}^{5}L_{6}^{i}.\nonumber
\end{eqnarray}
Similar to the estimates of $L_i,i=8,9,10$, we have
\begin{equation}\label{4.19}
\begin{split}
|L_{6}^{1}|&\leq C(1+t)^{-1}\int_{\mathbb{R}}(|\zeta|^{3}+|\phi|^{3})dx\\
&\leq C(1+t)^{-1}\left(\|\zeta(t)\|^{\frac{5}{2}}\|\zeta_{x}(t)\|^{\frac{1}{2}}+\|\phi(t)\|^{\frac{5}{2}}\|\phi_{x}(t)\|^{\frac{1}{2}}\right)
\\&\leq C\|(\zeta_{x},\phi_{x})(t)\|^2+C(1+t)^{-\frac{4}{3}},
\end{split}
\end{equation}
\begin{eqnarray}\label{4.20}
|L_{6}^{2}|+|L_{6}^{5}|&\leq& C\int_{\mathbb{R}}\left(|U_{x}|+|(\Theta_{-}^r)_t|+|(\Theta_{+}^r)_t|+|\Theta^c_t|+|\psi_{x}|\right)(|\zeta|^{3}+|\phi|^{3})dx
\nonumber\\&\leq& C\|(\zeta,\phi)(t)\|^{2}\|(\zeta_{x},\phi_{x})(t)\|\|(U^c_{x},\psi_{x})(t)\|+C\int_{\mathbb{R}}(\phi^2+\zeta^2)((U_{-}^{r})_x+(U_{+}^{r})_x)dx
\nonumber\\&\leq& C\left(\|(\zeta_{x},\phi_{x},\psi_{x})(t)\|^2+(1+t)^{-\frac{3}{2}}\right)+C\int_{\mathbb{R}}(\phi^2+\zeta^2)((U_{-}^{r})_x+(U_{+}^{r})_x)dx,
\end{eqnarray}
\begin{eqnarray}\label{4.21}
|L_{6}^{3}|&\leq& C\int_{\mathbb{R}}\left(|\zeta_{x}|+|\phi|(1+t)^{-\frac{1}{2}}\right)\left(\left|\left(\frac{g^{2}}{v}\right)_{x}\right|(\phi^{2}+\zeta^{2})+|(\phi(4\Gamma-\gamma P\phi))_{x}|\right)dx
\nonumber\\&\leq& C\int_{\mathbb{R}}\left(|\zeta_{x}|+|\phi|(1+t)^{-\frac{1}{2}}\right)\left(|\phi_{x}|+(1+t)^{-\frac{1}{2}}\right)(\phi^{2}+\zeta^{2})dx
\nonumber\\&&~~~~+C\int_{\mathbb{R}}\left(|\zeta_{x}|+|\phi|(1+t)^{-\frac{1}{2}}\right)(|\phi|+|\zeta|)(|\zeta_{x}|+|\phi_{x}|)dx
\nonumber\\&\leq& C\|\zeta_{x}(t)\|\left(\|\phi_{x}(t)\|+\|\zeta_{x}(t)\|+(1+t)^{-1}\right)+C(1+t)^{-1}\left(\|\phi_{x}(t)\|^{\frac{1}{2}}+\|\zeta_{x}(t)\|^{\frac{1}{2}}\right)
\nonumber\\&&~~~~+C(1+t)^{-\frac{1}{2}}\left(\|\phi_{x}(t)\|^{\frac{3}{2}}+\|\zeta_{x}(t)\|^{\frac{3}{2}}\right)
\nonumber\\&\leq& C(\|\phi_{x}(t)\|^{2}+\|\zeta_{x}(t)\|^{2})+C(1+t)^{-\frac{4}{3}},
\end{eqnarray}
\begin{equation}\label{4.22}
|L_{6}^{4}|\leq C\int_{\mathbb{R}}|\phi|(|\zeta|+|\phi|)(\omega^{2}+\omega_{x}^{2}+|R_2|)dx
\leq C\left(\|\omega(t)\|_1^{2}+\|(\phi_x,\zeta_x)(t)\|^2+(1+t)^{-\frac{7}{4}}\right).
\end{equation}
(\ref{4.5}) thus follows from (\ref{4.15})-(\ref{4.22}) and the smallness of  $\delta$ and $\eta$. The proof of Lemma 3.3 is completed.

For the estimate of $\|(\phi_{x},\psi_{x},\zeta_{x},\omega_{x})(t)\|$, we have
\begin{Lemma}Under the assumption Proposition 3.1, there exists a positive constant $C$ such that
\begin{equation}\label{3.30}
\begin{split}
&\|(\phi_{x},\psi_{x},\zeta_{x},\omega_{x})(t)\|^{2}+\int_{0}^{t}\|(\zeta_{xx},\omega_{xx},\omega_{x})(\tau)\|^2d\tau
\\&\leq C\left(\|(\phi_{0},\psi_{0},\zeta_{0},\omega_{0})\|_1^{2}+\delta^{\frac{1}{8}}+(\delta^{\frac{1}{8}}+\varepsilon+\eta)\int_{0}^{t}\|(\phi_{x},\psi_{x})(\tau)\|_1^{2}d\tau\right).
\end{split}
\end{equation}
\end{Lemma}
\noindent{\bf Proof.}~~Multiplying $(\ref{3.3})_{1x}$ by $\frac{P}{v}\phi_{x}$, $(\ref{3.3})_{2x}$ by $\psi_x$, $(\ref{3.3})_{3x}$ by $\frac{\zeta_{x}}{\theta}$, $(\ref{3.3})_{4x}$ by $\omega_x$, and adding the  resultant equations together, we have
\begin{equation}\label{3.31}
\begin{split}
&\left({\frac{P}{2v}\phi_{x}^{2}+\frac{\psi_{x}^{2}}{2}+\frac{R\zeta_{x}^{2}}{2(\gamma-1)\theta}+\frac{\omega_{x}^{2}}{2}}\right)_{t}+\frac{\kappa}{v\theta}\zeta_{xx}^{2}+\frac{A}{v}
\omega_{xx}^{2}+Av\omega_{x}^{2}+H_{2x}
\\&=R_{1x}\psi_{x}+R_{2x}\frac{\zeta_{x}}{\theta}+J_{1},
\end{split}
\end{equation}
where
\[
\begin{split}
H_{2}=&(p-P)_{x}\psi_{x}+\frac{\zeta_{x}}{\theta}\left((p-P)U_{x}-\kappa\left(\frac{\zeta_{x}}{v}-\frac{\Theta_{x}\phi}{vV}\right)_{x}\right)+\frac{\zeta_{x}w_{x}^{2}}{v\theta}
+\frac{v\omega^{2}\zeta_{x}}{\theta}-A\left(\frac{\omega_{x}}{v}\right)_{x}\omega_{x},\\
J_{1}=&\left(\frac{\zeta_{x}}{\theta}\right)_{x}\left(-\kappa\left(\frac{\zeta_{x}}{v}-\frac{\Theta_{x}\phi}{vV}\right)_{x}+(p-P)U_{x}\right)
+\frac{\kappa\zeta_{xx}^{2}}{v\theta}+\left(\frac{P}{2v}\right)_{t}\phi_{x}^{2}
+\frac{R\zeta_{x}^{2}}{2(\gamma-1)}\left(\frac{1}{\theta}\right)_{t}
\\&-p_{x}\psi_{x}\frac{\zeta_{x}}{\theta}-Av_x\omega\omega_{x}-A\omega_{xx}\omega_{x}\left(\frac{1}{v}\right)_{x}+\left(\frac{R}{v}\right)_{x}\zeta\psi_{xx}-\left(\frac{P}{v}\right)_{x}\phi\psi_{xx}.
\end{split}
\]
Integrating (\ref{3.31}) over $[0,t]\times\mathbb{R}$ yields
\begin{equation}\label{3.32}
\begin{split}
&\|(\phi_{x},\psi_{x},\zeta_{x},\omega_{x})(t)\|^{2}+\int_{0}^{t}\|(\zeta_{xx},\omega_{xx},\omega_{x})(\tau)\|^{2}d\tau
\\&\leq\|(\phi_{0x},\psi_{0x},\zeta_{0x},\omega_{0x})\|^{2}+\int_{0}^{t}\int_{R}J_{1}dxd\tau
+\int_{0}^{t}\int_{R}(|R_{1x}||\psi_{x}|+|R_{2x}||\zeta_{x}|)dxd\tau.
\end{split}
\end{equation}
It follows from  the  Cauchy inequality and (\ref{3.9})that
\begin{eqnarray}\label{3.33}
\int_{0}^{t}\int_{R}|J_{1}|dxd\tau&\leq&C(\eta+\varepsilon+\delta)\int_{0}^{t}\int_{\mathbb{R}}(\phi_x^2+\zeta_x^2+\zeta_{xx}^2+\psi_{xx}^2+\psi_x^2+w_x^2+w_{xx}^2)dxd\tau\nonumber\\
&&+C_\eta\underbrace{\int_{0}^{t}\int_{\mathbb{R}}|(\Theta_x,\Theta_{xx})|^2(\phi^2+\zeta^2)dxd\tau}_{J_1^1}.
\end{eqnarray}
Using (\ref{1.11}), (\ref{3.9}) and  Lemma 2.3, we obtain
\begin{eqnarray}\label{3.34}
J_1^1&\leq& \int_{0}^{t}\int_{\mathbb{R}}\left(\|((\Theta^r_\pm)_x,(\Theta^r_\pm)_{xx})(\tau)\|^2_{L^\infty}+|\Theta^c_x,\Theta^c_{xx})|^2\right)(\phi^2+\zeta^2)dxd\tau\nonumber\\
&\leq&\delta\int_{0}^{t}\int_{\mathbb{R}}(1+\tau)^{-1}e^{-\frac{c_{2}x^{2}}{1+\tau}}(\phi^{2}+\zeta^{2})dxd\tau
+\int_{0}^{t}\delta^\frac{1}{2}(1+\tau)^{-\frac{3}{2}}\sup_{0\leq\tau\leq t}\|(\phi,\zeta)(\tau)\|^2d\tau\nonumber\\
&\leq&\delta\int_{0}^{t}\int_{\mathbb{R}}(1+\tau)^{-1}e^{-\frac{c_{2}x^{2}}{1+\tau}}(\phi^{2}+\zeta^{2})dxd\tau
+C\varepsilon^2\delta^\frac{1}{2}.
\end{eqnarray}
Similar to the estimate of (\ref{3.28}), we have
\begin{equation}\label{3.35}
  \begin{split}
     \int_0^t\int_{\mathbb R}\left(|R_{1x}\psi_x|+|R_{2x}\zeta_x|\right)dxd\tau
     \le &\delta^{\frac{1}{8}}\int_0^t\|(\psi_{xx},\zeta_{xx})(\tau)\|^2d\tau+C\delta^{\frac{1}{8}}.
  \end{split}
\end{equation}
Combining (\ref{3.32})-(\ref{3.35}) and using Lemmas 3.2-3.3 and the smallness of $\eta,\varepsilon$ and $\delta$ leads to (\ref{3.30}) immediately. The proof of Lemma 3.4 is finished.

The next lemma give the estimate on $\|(\phi_{xx},\psi_{xx},\zeta_{xx},\omega_{xx})(t)\|$.
\begin{Lemma}Under the assumption Proposition 3.1, there exists a positive constant $C$ such that
\begin{equation}\label{3.36}
\begin{split}
&\|(\phi_{xx},\psi_{xx},\zeta_{xx},\omega_{xx})(t)\|^{2}+\int_{0}^{t}\|(\zeta_{xxx},\omega_{xxx},\omega_{xx})(\tau)\|^2d\tau
\\&\leq C\left(\|(\phi_{0},\psi_{0},\zeta_{0},\omega_{0})\|_2^{2}+\delta^{\frac{1}{8}}+(\delta^{\frac{1}{8}}+\varepsilon+\eta)\int_{0}^{t}\|(\phi_{x},\psi_{x})(\tau)\|_1^{2}d\tau\right).
\end{split}
\end{equation}
\end{Lemma}
\noindent{\bf Proof.}~~Multiplying $(\ref{3.3})_{1xx}$ by $\frac{P}{v}\phi_{xx}$, $(\ref{3.3})_{2xx}$ by $\psi_{xx}$, $(\ref{3.3})_{3xx}$ by $\frac{\zeta_{xx}}{\theta}$, $(\ref{3.3})_{4xx}$ by $\omega_{xx}$, and adding the  resultant equations together, we obtain
\begin{equation}\label{3.37}
\begin{split}
&\left({\frac{P}{2v}\phi_{xx}^{2}+\frac{\psi_{xx}^{2}}{2}+\frac{R\zeta_{xx}^{2}}{2(\gamma-1)\theta}+\frac{\omega_{xx}^{2}}{2}}\right)_{t}+\frac{\kappa}{v\theta}\zeta_{xxx}^{2}+\frac{A}{v}
\omega_{xxx}^{2}+Av\omega_{xx}^{2}+H_{3x}
\\&=J_{2}+J_{3}+J_4+J_4+J_5+R_{1xx}\psi_{xx}+R_{2xx}\frac{\zeta_{xx}}{\theta},
\end{split}
\end{equation}
where
\[
\begin{split}
&H_{3}=(p-P)_{xx}\psi_{xx}+\frac{\zeta_{xx}}{\theta}\left((p-P)U_{x}-\kappa\left(\frac{\zeta_{x}}{v}-\frac{\Theta_{x}\phi}{vV}\right)_{x}\right)_{x}
-A\left(\frac{w_{x}}{v}\right)_{xx}\omega_{xx},\\
&J_{2}=\left(\frac{P}{2v}\right)_{t}\phi_{xx}^{2}+\left(\frac{R}{2(\gamma-1)\theta}\right)_{t}\zeta_{xx}^{2}-(2p_{x}\psi_{xx}+p_{xx}\psi_{x})\frac{\zeta_{xx}}{\theta},\\
&J_{3}=\frac{\kappa\theta_x\zeta_{xx}\zeta_{xxx}}{\theta^2v}+\kappa\left(\frac{\theta_x\zeta_{x}}{\theta^2}
-\frac{\zeta_{xxx}}{\theta}\right)\left[\frac{-2v_x\zeta_{xx}}{v^2}+\zeta_x\left(\frac{-v_x}{v^2}\right)_x-\left(\frac{\Theta_x\phi}{vV}\right)_{xx}-((p-P)U_x)_x\right],\\
&J_{4}=\frac{\zeta_{xx}}{\theta}\left(\frac{w_x^2}{v}+vw^2\right)_{xx}+
\frac{2Aw_{xxx}w_{xx}v_x}{v^2}+Aw_{xxx}w_{x}\left(\frac{v_x}{v^2}\right)_x-Aww_{xx}v_{xx}-2Aw_{xx}w_xv_{x},\\
&J_{5}=\left(\frac{R}{v}\right)_{xx}\zeta\psi_{xxx}-\left(\frac{P}{v}\right)_{xx}\phi\psi_{xxx}
+2\left(\frac{R}{V}\right)_{x}\zeta_x\psi_{xxx}-2\left(\frac{P}{v}\right)_{x}\phi_x\psi_{xxx}.
\end{split}
\]
The Cauchy inequality,(\ref{1.1}), (\ref{1.11}), (\ref{3.9}) and  Lemma 2.3 imply that
\begin{eqnarray}\label{3.38}
\int_0^t\int_{\mathbb{R}}\sum_{i=2}^4|J_i|dxd\tau&\leq& C(\varepsilon+\delta)\int_{0}^{t}\left(\|(\zeta_x,w_x)(\tau)\|_2^2+\|(\phi_x,\psi_x)(\tau)\|_1^2\right)d\tau\nonumber\\
&&+C\int_{0}^{t}\int_{\mathbb{R}}(\Theta_{x}^2+\Theta_{xx}^2)(\phi^2+\zeta^2)dxd\tau+\eta\int_{0}^{t}\|(w_{xx},\zeta_x)(\tau)\|^2d\tau\nonumber\\
&&+C\underbrace{\int_{0}^{t}\int_{\mathbb{R}}\left(|\zeta_{xx}|^3+|\zeta_{xx}w_{xx}^2|\right)dxd\tau}_{K_1}\nonumber\\
&&+C_\eta\underbrace{\int_{0}^{t}\int_{\mathbb{R}}\left(|\Theta_{x}|^8+|\Theta_{x}^2\Theta_{xxx}^2|+|\Theta_{x}^4\Theta_{xx}^2|+w^2V_{xx}^2\right)dxd\tau}_{K_2}.
\end{eqnarray}
The second term on the right hand side of (\ref{3.38}) has been estimated in (\ref{3.34}). For $K_1$ and $K_2$, we deduce from the Cauchy inequality, the Sobolev inequality,  (\ref{1.11}), (\ref{3.9}) and  Lemma 2.3 that
\begin{equation}\label{3.39}
\begin{split}
K_1\leq& \int_{0}^{t}\sup_{0\leq\tau\leq t}\{\|\zeta_{xx}(\tau)\|\}\|\zeta_{xx}(\tau)\|^{\frac{3}{2}}\|\zeta_{xxx}(\tau)\|^{\frac{1}{2}}d\tau\\
&+\int_{0}^{t}\sup_{0\leq\tau\leq t}\{\|w_{xx}(\tau)\|^2\}\|w_{xx}(\tau)\|^{\frac{1}{2}}\|w_{xxx}(\tau)\|^{\frac{1}{2}}\|\zeta_{xx}(\tau)\|d\tau\\
\leq&C\varepsilon\int_{0}^{t}\|(\zeta_{xx}, w_{xx})\|_1^2d\tau,
\end{split}
\end{equation}
\begin{equation}\label{3.40}
\begin{split}
K_2\leq& \int_{0}^{t}\left(\|((V_+^r)_{xx},(V_-^r)_{xx})(\tau)\|_{L^\infty}^2+\|\Theta_{xx}^c(\tau)\|^2_{L^\infty}\right)\sup_{0\leq\tau\leq t}\{\|w(\tau)\|^2\}d\tau\\
&+\int_{0}^{t}\int_{\mathbb{R}}\left((\Theta_{x}^c)^8+(\Theta_+^r)_{x})^4+(\Theta_-^r)_{x})^4\right)dxd\tau\\
\leq&C\int_{0}^{t}\left(\varepsilon^2\delta^{\frac{1}{2}}(1+\tau)^{-\frac{3}{2}}+\delta^2(1+\tau)^{-\frac{5}{2}}+\delta(1+\tau)^{-3}\right)d\tau\\
\leq&C\delta^{\frac{1}{2}}.
\end{split}
\end{equation}

To estimate $J_5$, we have from integration by parts and $(\ref{3.3})_1$ that
\begin{eqnarray}\label{3.41}
\int_{0}^{t}\int_{\mathbb{R}}\left(\frac{P}{v}\right)_{xx}\phi\psi_{xxx}dxd\tau
&=&\int_{0}^{t}\int_{\mathbb{R}}\left\{\frac{-P\phi}{2v^2}\phi^2_{xx}\right\}_tdxd\tau+
\int_{0}^{t}\int_{\mathbb{R}}\left\{\frac{-P\phi}{2v^2}\right\}_t\phi^2_{xx}dxd\tau\nonumber\\ &&-\int_{0}^{t}\int_{\mathbb{R}}\psi_{xx}\left\{\phi\left[\frac{P_{xx}}{v}+2P_{x}\left(\frac{1}{v}\right)_{x}+P\left(\frac{2v_{x}^{2}}{v^{3}}-\frac{V_{xx}}{v^{2}}\right)\right]\right\}_xdxd\tau\nonumber\\
&\leq&-\int_{\mathbb{R}}\frac{P\phi}{2v^2}\phi_{xx}^2dx+\int_{\mathbb{R}}\frac{P(V_0,\Theta_0)\phi_0}{2v_0^2}\phi_{0xx}^2dx+\eta\int_{0}^{t}\left\|\psi_{xx}(\tau)\right\|^2d\tau\nonumber\\
&&+C(\delta+\varepsilon)\int_{0}^{t}\int_{\mathbb{R}}\left(\psi_{xx}^2+\phi_{x}^2+\phi_{xx}^2\right)dxd\tau+C_\eta\int_{0}^{t}\int_{\mathbb{R}}\Theta_{x}^2\phi^2dxd\tau\nonumber\\
&\leq&C\left(\varepsilon\|\phi_{xx}(t)\|^2+\|\phi_{0xx}\|^2\right)+C(\delta+\varepsilon)\int_{0}^{t}\left\|(\psi_{xx},\phi_{x},\phi_{xx})(\tau)\right\|^2d\tau\nonumber\\
&&+\eta\int_{0}^{t}\left\|\psi_{xx}(\tau)\right\|^2d\tau+C_\eta\int_{0}^{t}\int_{\mathbb{R}}\Theta_{x}^2\phi^2dxd\tau.
\end{eqnarray}
Similarly, it holds that
\begin{equation}\label{3.42}
\begin{split}
&\int_{0}^{t}\int_{\mathbb{R}}\left(\frac{R}{v}\right)_{xx}\zeta\psi_{xxx}dxd\tau\\
\leq &C\left(\varepsilon\|\phi_{xx}(t)\|^2+\|\phi_{0xx}\|^2\right)+\eta\int_{0}^{t}\left\|(\zeta_{xx},\zeta_{xxx},\psi_{xx})(\tau)\right\|^2d\tau\\
&+C(\delta+\varepsilon)\int_{0}^{t}\left\|(\psi_{xx},\zeta_x,\phi_{x})(\tau)\right\|^2d\tau+C_\eta\int_{0}^{t}\int_{\mathbb{R}}(\varepsilon^2\phi_{xx}^2+\Theta_x^2\zeta^2)dxd\tau,
\end{split}
\end{equation}
and \begin{equation}\label{3.43}
\begin{split}
\int_{0}^{t}\int_{\mathbb{R}}\left(\left(\frac{R}{V}\right)_{x}\zeta_x\psi_{xxx}-\left(\frac{P}{v}\right)_{x}\phi_x\psi_{xxx}\right)dxd\tau
\leq
C(\delta+\varepsilon)\int_{0}^{t}\left\|(\phi_{x},\zeta_{x}, \zeta_{xx}, \psi_{xx},\phi_{xx})(\tau)\right\|^2d\tau.
\end{split}
\end{equation}
Therefore, it follows from (\ref{3.41})-(\ref{3.43}) that
\begin{equation}\label{3.44}
\begin{split}
\int_{0}^{t}\int_{\mathbb{R}}|J_5|dxd\tau
\leq&
C(\delta+\varepsilon)\int_{0}^{t}\left\|(\phi_{x},\zeta_{x}, \zeta_{xx}, \psi_{xx},\phi_{xx})(\tau)\right\|^2d\tau+\eta\int_{0}^{t}\left\|(\zeta_{xx},\zeta_{xxx},\psi_{xx})(\tau)\right\|^2d\tau\\
&+C\left(\varepsilon\|\phi_{xx}(t)\|^2+\|\phi_{0xx}\|^2\right)+C_\eta\int_{0}^{t}\int_{\mathbb{R}}\left(\varepsilon^2\phi_{xx}^2+\Theta_x^2(\zeta^2+\phi^2)\right)dxd\tau.
\end{split}
\end{equation}
Finally, by using Lemma 3.1, the Sobolev inequality and the Young inequality,  we have
\begin{eqnarray}\label{3.45}
&&\int_{0}^{t}\int_{\mathbb{R}}\left|R_{1xx}\psi_{xx}+R_{2xx}\frac{\zeta_{xx}}{\theta}\right|dxd\tau\nonumber\\
&\leq&
C\int_{0}^{t}\left(\|\psi_{xx}(\tau)\|\|R_{1xx}(\tau)\|+\|\zeta_{xx}(\tau)\|^{\frac{1}{2}}\|\zeta_{xxx}(\tau)\|^{\frac{1}{2}}\|R_{2xx}(\tau)\|_{L^1}\right)d\tau\nonumber\\
&\leq&
C\int_{0}^{t}\left(\|\psi_{xx}(\tau)\|\delta(1+\tau)^{-\frac{5}{4}}+\|\zeta_{xx}(\tau)\|^{\frac{1}{2}}\|\zeta_{xxx}(\tau)\|^{\frac{1}{2}}\delta^{\frac{1}{8}}(1+\tau)^{-\frac{7}{8}}\right)d\tau\nonumber\\
&\leq&
C\delta^{\frac{1}{8}}\int_{0}^{t}\|(\psi_{xx},\zeta_{xx},\zeta_{xxx})(\tau)\|^2d\tau+C\delta^{\frac{1}{8}}.
\end{eqnarray}

Integrating (\ref{3.37}) over $[0,t]\times\mathbb{R}$, and using (\ref{3.38})-(\ref{3.40}), (\ref{3.44})-(\ref{3.45}), Lemmas 3.2-3.4 and the smallness of $\eta, \varepsilon$ and $\delta$, we  get (\ref{3.36}) immediately. This completes the proof of Lemma 3.6.

As a direct consequence of Lemmas 3.2-3.5, we have the following corollary.
\begin{Corollary}Under the assumption Proposition 3.1, there exists a positive constant $C$ such that
\begin{equation}\label{3.46}
\begin{split}
&\|(\phi,\psi,\omega,\zeta)(t)\|_{2}^{2}+\int_{0}^{t}\left(\|\zeta_x(\tau)\|_{2}^{2}+\|\omega(\tau)\|_{3}^{2}\right)d\tau
\\
\leq &C\left(\|(\phi_{0},\varphi_{0},\omega_{0},\zeta_{0})\|_{2}^{2}+\delta^{\frac{1}{8}}+
(\delta^{\frac{1}{8}}+\varepsilon+\eta)\int_{0}^{t}\|(\phi_{x},\psi_{x})(\tau)\|_{1}^{2}d\tau\right).
\end{split}
\end{equation}
\end{Corollary}

For the reminder term $\displaystyle\int_{0}^{t}\|(\phi_{x},\psi_{x})(\tau)\|_1^2d\tau$ in (\ref{3.46}), we have
\begin{Lemma}Under the assumption Proposition 3.1, there exists a positive constant $C$ such that
\begin{equation}\label{3.47}
\int_{0}^{t}\|(\phi_{x},\psi_{x})(\tau)\|_1^2d\tau\leq C\left(\|(\phi_{0},\psi_{0},\zeta_{0},\omega_{0})\|_2^{2}+\delta^{\frac{1}{8}}\right).
\end{equation}
\end{Lemma}
\noindent{\bf Proof.}~~Multiplying $(\ref{3.3})_{2}$ by $-\frac{P}{2}\phi_{x}$, $(\ref{3.3})_{3}$ by $\psi_{x}$, and adding the  resultant equations together, we have
\begin{equation}\label{3.48}
\begin{split}
&\displaystyle\left\{\frac{R}{\gamma-1}\zeta\psi_{x}-\frac{P}{2}\phi_{x}\psi\right\}_{t}+\left\{\frac{P}{2}\phi_{t}\psi-\frac{R}{\gamma-1}\zeta\psi_{t}\right\}_{x}
+\frac{P^{2}}{2v}\phi_{x}^{2}+\frac{P}{2}\psi_{x}^{2}
\\=&\frac{P_{x}}{2}\psi\psi_{x}-\frac{P_{t}}{2}\phi_{x}\psi
+\frac{P}{2}\phi_{x}\left(\frac{R\zeta_x}{v}-\frac{(R\zeta-P\phi)v_x}{v^2}-R_{1}\right)-\frac{R}{\gamma-1}\zeta_{x}\psi_{t}
\\&+\kappa\left(\frac{\zeta_{x}}{v}-\frac{\Theta_{x}\phi}{vV}\right)_{x}\psi_{x}-(p-P)(U_{x}+\psi_{x})\psi_{x}
+R_{2}\psi_{x}+\frac{\omega_{x}^{2}}{v}\psi_{x}+v\omega^{2}\psi_{x}.
\end{split}
\end{equation}
Integrating (\ref{3.48}) on $[0,t]\times\mathbb{R}$ yields
\begin{equation}\label{3.49}
\begin{split}
&\int_{0}^{t}\int_{\mathbb{R}}(\phi_{x}^{2}+\psi_{x}^{2})dxd\tau
\\\leq&C\|(\phi,\psi,\zeta,\omega)(t)\|_{1}^{2}+C\|(\phi_{0},\psi_{0},\zeta_{0},\omega_{0})\|_{1}^{2}
+C\int_{0}^{t}\|(\zeta_{x},\omega_{x})(\tau)\|_1^{2}d\tau
\\&+(\frac{1}{4}+\delta+\varepsilon)\int_{0}^{t}\int_{\mathbb{R}}(\phi_{x}^{2}+\psi_{x}^{2})dxd\tau+
C\int_{0}^{t}\int_{\mathbb{R}}(|\Theta_{x}|+|\Theta_{xx}|)^{2}(\phi^{2}+\psi^{2})
\\&+C\int_{0}^{t}\int_{R}|R_{1}\phi_{x}+R_{2}\psi_{x}|dxd\tau.
\end{split}
\end{equation}
Similar to (\ref{3.45}), we have
\begin{eqnarray}\label{3.50}
\int_{0}^{t}\int_{\mathbb{R}}\left|R_{1}\phi_{x}+R_{2}\psi_{x}\right|dxd\tau
&\leq&
C\int_{0}^{t}\left(\|\phi_{x}(\tau)\|\|R_{1}(\tau)\|+\|\psi_{x}\|^{\frac{1}{2}}\|\psi_{xx}\|^{\frac{1}{2}}\|R_{2}(\tau)\|_{L^1}\right)d\tau\nonumber\\
&\leq&
C\int_{0}^{t}\left(\|\phi_{x}(\tau)\|\delta(1+\tau)^{-\frac{5}{4}}+\|\psi_{x}(\tau)\|^{\frac{1}{2}}\|\psi_{xx}(\tau)\|^{\frac{1}{2}}\delta^{\frac{1}{8}}(1+\tau)^{-\frac{7}{8}}\right)d\tau\nonumber\\
&\leq&
C\delta^{\frac{1}{8}}\int_{0}^{t}\|(\phi_{x},\psi_{x},\psi_{xx})(\tau)\|^2d\tau+C\delta^{\frac{1}{8}}.
\end{eqnarray}
Combining (\ref{3.49})-(\ref{3.50}), and  using (\ref{3.29}), (\ref{3.34}) and  (\ref{3.46}), we get by the smallness of $\varepsilon,\delta$ and $\eta$ that
\begin{equation}\label{3.51}
\int_{0}^{t}\int_{\mathbb{R}}(\phi_{x}^{2}+\psi_{x}^{2})dxd\tau
\leq C\left(\|(\phi_{0},\psi_{0},\zeta_{0},\omega_{0})\|_{2}^{2}+\delta^{\frac{1}{8}}\right)+C(\varepsilon+\delta^{\frac{1}{8}})\int_{0}^{t}\|(\phi_{xx},\psi_{xx})(\tau)\|^2d\tau.
\end{equation}

Similarly, by multiplying $(\ref{3.3})_{2x}$ by $-\frac{P}{2}\phi_{xx}$, $(\ref{3.3})_{3x}$ by $\psi_{xx}$, and repeating the same argument as above, we can also obtain
\begin{equation}\label{3.52}
\int_{0}^{t}\int_{\mathbb{R}}(\phi_{xx}^{2}+\psi_{xx}^{2})dxd\tau
\leq C\left(\|(\phi_{0},\psi_{0},\zeta_{0},\omega_{0})\|_{2}^{2}+\delta^{\frac{1}{8}}\right)+C(\varepsilon+\delta^{\frac{1}{8}})\int_{0}^{t}\|(\phi_{x},\psi_{x})(\tau)\|^2d\tau.
\end{equation}
(\ref{3.47}) thus follows from (\ref{3.51})-(\ref{3.52}) and the smallness of $\varepsilon$ and $\delta$. This completes the proof of Lemma 3.6.
\\

\noindent{\bf Proof of Proposition 3.1.}~~Proposition 3.1 follows from Corollary 3.1, Lemma 3.6, and the smallness of $\varepsilon,\delta$ and $\eta$ immediately.

\begin{center}
{\bf Acknowledgement}
\end{center}
This work
was supported by the National Natural Science Foundation of China
(Grant No. 11501003),  and the Cultivation Fund
of Young Key Teacher at Anhui University.

\end{document}